\documentclass[12pt,leqno]{article}
\usepackage{amssymb,amsmath,amsthm}
\usepackage[all]{xy}
\usepackage[backend=bibtex,style=alphabetic,sorting=nyt]{biblatex}
\numberwithin{equation}{section}
\usepackage[top=2cm,bottom=2cm,left=2cm,right=4cm,marginparsep=0.3cm,marginparwidth=3cm,includefoot]{geometry}
\usepackage{mathtools}
\usepackage{radon}

\usepackage{enumerate}
\usepackage[colorlinks=true, pdfstartview=FitV, linkcolor=blue, citecolor=blue, urlcolor=blue]{hyperref}
\usepackage[normalem]{ulem}  
\usepackage{mathrsfs}
\usepackage{yfonts}
\usepackage{accents}
\usepackage{stmaryrd}
\usepackage{xspace}
\usepackage{marginnote}
\usepackage{tikz-cd}
\usepackage{tocvsec2}
\usepackage{rotating}
\usepackage{parallel} 
\usepackage{framed}
\usepackage{lipsum}

\numberwithin{equation}{section}
\usepackage[framemethod=default]{mdframed}
\usepackage{tikz}
\usetikzlibrary{patterns}
\usepackage{graphicx}
\usepackage{stmaryrd}
\usepackage{geometry}

\newtheorem*{theorem*}{Theorem}

\newtheorem{proposition-definition}{Proposition-definition}

\makeatletter
\renewcommand{\l@section}{\@dottedtocline{1}{1.5em}{2.6em}}
\renewcommand{\l@subsection}{\@dottedtocline{2}{4.0em}{3.6em}}
\renewcommand{\l@subsubsection}{\@dottedtocline{3}{7.4em}{4.5em}}
\makeatother

\UseRawInputEncoding

\DeclareMathAlphabet{\mathpzc}{OT1}{pzc}{m}{it}
\DeclareSymbolFontAlphabet{\amsmathbb}{AMSb}%

\newmdenv[
frametitle=Reminder,
skipabove=\topsep,
skipbelow=\topsep,
]{reminder}

\newmdenv[
frametitle=Problem,
skipabove=\topsep,
skipbelow=\topsep,
]{Problem}

\newcommand\restr[2]{{
	\left.\kern-\nulldelimiterspace
	#1
	\right|_{#2} 
}}

\begin{document}
	\title
	{Action of micro-differential operators on quantized contact transformations}
	\author{Mehdi Benchoufi}
	\maketitle

\begin{abstract}
	Quantized contact transforms (QCT) have been constructed in~\cite{SKK73}. We give here a complete proof of the fact that such QCT commute with the action of microdifferential operators. To our knowledge, such a proof did not exist in the literature. We apply this result to the microlocal Radon transform.
\end{abstract}

\tableofcontents

\section{Introduction}
\subsection{Overview of the results}
	
For a manifold $M$, let us denote by $T^*M$ the cotangent bundle and $\dT{}^*M$ the bundle $T^*M$ with the zero section removed. We will consider the following situation: let $X$ and $Y$ be two complex manifolds of the same dimension, a closed submanifold $Z$ of $X\times Y$, open subset $U\subset \dT{}^*X$ and $V\subset \dT{}^*Y$ and assume that the conormal bundle $\dT{}^*_Z(X\times Y)$  induces a contact transformation 
\eq\label{eq:diag2}
&&\xymatrix{
	&\dT{}^*_Z(X\times Y)\cap (U\times V^a)\ar[ld]_-{\sim}\ar[rd]^-{\sim}&\\
	\dT{}^*X\supset U\ar[rr]^-{\sim}&&V\subset \dT{}^*Y.
}
\eneq
Let $F\in\Derb(X)$ and let $\phi_K(F)$ denote the contact transformation of $F$ with kernel $K\in\Derb(X\times Y)$. We will prove an isomorphism between $\muhom(F,\sho_X)$ on $U\cap T_M^*X$ and $\muhom(\phi_K(F),\sho_Y)$ on $V\cap T_N^*Y$, which follows immediately from \cite[Lem. 11.4.3]{KS90}. Our main result will be the commutation of this isomorphism to the action of microdifferential operators. Although considered as well-known, the proof of this  commutation does not appear clearly in the literature
(see \cite[p.~467]{SKK73}), and is far from being obvious. In fact, we will consider a more general setting, replacing sheaves of microfunctions with sheaves of the type $\muhom(F,\sho_X)$. Moreover, assume that $X$ and $Y$ are complexification of real analytic manifolds $M$ and $N$ respectively.  Then, it is known that, under suitable hypotheses, one can quantize this contact transform  and get an isomorphism between microfunctions on $U\cap T_M^*X$ and microfunctions on $V\cap T_N^*Y$ \cite{KKK86}. 

Then, we will specialize our results to the case of projective duality. We will study the {\em microlocal} Radon transform understood as a quantization of projective duality, both in the real and the complex case. 

In the real case, denote by $P$ the real projective space (say of dimension $n$), by $P^*$ its dual and by $S$ the incidence relation:
\eq\label{eq:incidence1}
&&S\eqdot\{(x,\xi)\in P\times P^*; \langle x,\xi\rangle=0\}.
\eneq
In this setting, there is a well-known correspondance between distributions on $P$ and $P^*$ due to  Gelfand, Gindikin, Graev \cite{GGG}\, and to Helgason \cite{Hel}. However, it is known since the 70s under the influence of the Sato's school, that to well-understand what happens on real (analytic) manifolds, it may be worth to look at their complexification. 

Hence,  denote by $\BBP$ the complex projective space of dimension $n$,  by $\BBP^*$ the dual projective space and by $\BBS\subset\BBP\times\BBP^*$ the incidence relation. We have the correspondence
\eq\label{eq:diag1}
&&\xymatrix@C-0pc@R+0pc{
	&\dT{}^*_{\mathbb{S}}(\mathbb{P}\times\mathbb{P}^*)\ar[ld]_-{\sim}\ar[rd]^-{\sim}&\\
	\dT{}^*\mathbb{P}\ar[rr]^-{\sim}&&\dT{}^*\mathbb{P}^*
}
\eneq

This contact transformation induces an equivalence of categories between perverse sheaves modulo constant ones on the complex projective space and perverse sheaves modulo constant ones on its dual, as shown by Brylinski \cite{B86}, or between coherent $\mathcal{D}$-modules modulo flat connections, as shown by D'Agnolo-Schapira \cite{DS94}.

In continuation of the previous cited works, we shall consider the contact transform induced by~\eqref{eq:diag1}
\eq\label{eq:diag2}
&&\xymatrix{
	&\dT{}^*_\BBS(\BBP\times\BBP^*)\cap (\dT{}^*_P\BBP\times \dT{}^*_{P^*}\BBP^*)\ar[ld]_-{\sim}\ar[rd]^-{\sim}&\\
	\dT{}_P^*\BBP\ar[rr]^-{\sim}&&\dT{}_{P^*}^*\BBP^*
}
\eneq
The above contact transformation leads to the well-known fact that the Radon transform establishes an isomorphism of sheaves of microfunctions on $P$ and $P^*$ (see~\cite{KKK86}). We will apply our main result to prove the commutation of this isomorphism to the action of microdifferential operators.

\subsection{Main theorems}\label{The problem}

We will use the langage of sheaves and $\mathcal{D}$-modules and we refer the reader to \cite{KS90} and \cite{K03} for a detailed developement of these topics. We will denote by k a commutative unital ring of finite global dimension.

\paragraph*{Notations for integral transforms}\label{def:radon_transforms}
Let $X,Y$ be real manifolds and $S$ a closed submanifold $X\times Y$. Consider the diagrams $X \from[f] S \to[g] Y$, $X \from[q_1] X\times Y \to[q_2] Y$. Let $F\in\Derb(\cor_Y)$ and $K\in\Derb(\cor_{X\times Y})$. The integral transform of $F$ with respect to the kernel $K$ is defined to be $\Phi_{K}(F):=\reim{q_1}(K\tens\opb{q_2}F)$. We will denote $\Phi_{S}(F)$ the integral transform of $F$ with respect to the kernel $\cor_S[d_S-d_X]$.

\paragraph*{Results on the functor $\muhom$ }\label{sec:pre_notations_results}

To establish our main results, we will need the following complement on the functor $\muhom$.

For $(M_i)_{i=1,2,3}$, three manifolds, we write $M_{ij}\eqdot M_i\times M_j$ ($1\leq i,j\leq3$). We consider the operation of composition of kernels:
\eq\label{eq:conv}
&&\ba{l}
\conv[2]\;\cl\;\Derb(\cor_{M_{12}})\times\Derb(\cor_{M_{23}})\to\Derb(\cor_{M_{13}})\\
\hs{10ex}\ba{rcl}(K_1,K_2)\mapsto K_1\conv[2] K_2&\eqdot&
\reim{q_{13}}(\opb{q_{12}}K_1\tens\opb{q_{23}}K_2)\\
&\simeq&\reim{q_{13}}\opb{\delta_2}(K_1\etens K_2).\ea
\ea
\eneq

We add a subscript $a$ to $p_j$ to denote by $p_j^a$ the composition of $p_j$ and the antipodal map on $T^*M_{j}$. We define the composition of kernels on cotangent bundles
(see~\cite[Prop.~4.4.11]{KS90})
\eq\label{eq:aconv}
&&\hs{-0ex}\ba{rcl}
\aconv[2]\;\cl\;\Derb(\cor_{T^*M_{12}})\times\Derb(\cor_{T^*M_{23}})
&\to&\Derb(\cor_{T^*M_{13}})\\
(K_1,K_2)&\mapsto&K_1\aconv[2] K_2\eqdot
\reim{p_{13}}(\opb{p_{12^a}} K_1\tens\opb{p_{23}} K_2)\\
\ea
\eneq
	Let $F_i,G_i,H_i$ respectively in $\Derb(\cor_{M_{12}}),\Derb(\cor_{M_{23}}),\Derb(\cor_{M_{34}})$, $i=1,2$.
	Let $U_i$ be an open subset of $T^*M_{ij}$ $\lp i=1,2$, $j=i+1\rp$ and set 
	\eqn
	U_3=U_i\aconv[2]U_j=p_{13}(\opb{p_{12^a}}(U_1)\cap\opb{p_{23}}(U_2))
	\eneqn 
	In \cite{KS90}, a canonical morphism in $\Derb(\cor_{T^*M_{13}})$ is constructed
	\eq\label{aconv_natural_morphism}
	&&\muhom(F_1,F_2)\vert_{U_1}\aconv[2]\muhom(G_1,G_2)\vert_{U_2}\to\muhom(F_1\conv[2]G_1,F_2\conv[2]G_2)\vert_{U_3}.
	\eneq
	We will see that the composition $\aconv$ is associative and we will see also that the morphism (\ref{aconv_natural_morphism}) is compatible with associativity with respect to $\aconv$.
  

\paragraph*{Complex contact transformations}
Consider now two complex manifolds $X$ and $Y$ of the same dimension $n$, open $\mathbb{C}^{\times}$-conic subsets $U$ and $V$ of $\dT X$ and $\dT Y$, respectively, $\Lambda$ a smooth closed submanifold of $U\times V^a$ and assume that the projections $p_1\vert_\Lambda$ and $p_2^a\vert_\Lambda$ induce isomorphisms, hence a homogeneous symplectic isomorphism $\chi\cl U\isoto V$:

\eqn
&&\xymatrix{
	&\Lambda\subset U\times V^a\ar[ld]^-{p_1}_-{\sim}\ar[rd]_-{p_2^a}^-{\sim}&\\
	\dTX\supset U\ar[rr]^-{\sim}_-{\chi}&  &V\subset \dTY
}
\eneqn

Let us consider a perverse sheaf $L$ on $X\times Y$ satisfying $(\opb{p_1}(U)\cup\opb{{p_2^a}}(V))\cap\SSi(L)\subset\Lambda$ and a section $s$ of $\muhom(L,\Omega_{X\times Y/X})$ on $\Lambda$, where $\Omega_{X\times Y/X}:=\mathcal{O}_{X\times Y}\tens_{\opb{q_2}\mathcal{O}_{Y}}\opb{q_2}\Omega_Y$. Recall that one denotes by $\she_X^{\mathbb{R}}$ sheaf of rings $\she_X^{\mathbb{R}}:=\muhom(\mathbb{C}_{\Delta_{X}},\Omega_{X\times X/X})[d_X]$, and $\she_{X}$ the subsheaf of $\she_X^{\mathbb{R}}$ of finite order microdifferential operators. In the following theorem, the first statement $(i)$ is well-known, see \cite{SKK73}, $(ii)$ is proved in \cite{KS90}, the fact that isomorphism (\ref{eq:intro_qct_main_theorem}) is compatible with the action of microdifferential operators was done at the germ level in \cite{KS90}, but from a global perspective, it was announced for microfunctions in various papers but no detailed proof exists  to our knowledge. We will prove our main theorem:

\begin{theorem}\label{th:microfunction_contact_iso}
	Let $G\in\Derb(\C_Y)$ and assume to be given a section $s$ of $\muhom(L,\Omega_{X\times Y/X})$, non-degenerate on $\Lambda$. 
	\bnum
	\item For $W\subset U$, $P\in\she_X(W)$, there is a unique $Q\in\she_Y(\chi(W))$ satisfying $P\cdot s=s\cdot Q$ ($P,Q$ considered as sections of $\she_{X\times Y}$). The morphism induced by $s$
	\eqn
	\opb{\chi}\she_Y\vert_{V}\to\she_X\vert_{U}
	\eneqn
	\vspace{-2.5em}
	\eqn
	P\mapsto Q
	\eneqn
	is a ring isomorphism. 
	\item We have the following isomorphism in $\Derb(\C_U)$
	\eq\label{eq:intro_qct_main_theorem}
	&&\opb{\chi}\muhom(G,\sho_Y)\vert_V\isoto\muhom(\Phi_{L[n]}(G),\sho_X)\vert_U
	\eneq
	\item The isomorphism (\ref{eq:intro_qct_main_theorem}) is compatible with the action of $\she_Y$ and $\she_X$ on the left and right side of (\ref{eq:intro_qct_main_theorem}) respectively.
	\enum
\end{theorem}

We will see that the action of microdifferential operators in Theorem \ref{th:microfunction_contact_iso} (iii) is derived from the morphism (\ref{aconv_natural_morphism}).

\paragraph*{Projective duality for microfunctions}\label{sec:notation_from_results}

For $M$ a real analytic manifold and $X$ its complexification, we might be led to identify $T_M^*X$ with $i\cdot T^*M$. We denote by $\sha_M$, $\shb_M$, $\shc_M$ the sheaves of real analytic functions, hyperfunctions, microfunctions, respectively.

In this article, we will quantize the contact transform associated with the Lagrangian submanifold $\dT{}_\mathbb{S}^*(\mathbb{P}\times\mathbb{P}^*)$. We will construct and denote by $\chi$ the homogeneous symplectic isomorphism between $\dT{}^*\mathbb{P}$ and $\dT{}^*\mathbb{P}^*$. 

For $\varepsilon\in\mathbb{Z}/2\mathbb{Z}$, we denote by $\mathbb{C}_{P}(\varepsilon)$ the two locally constant sheaf of rank one on $P$ (see Section \ref{sec:notations_projective} for a precise definition).

Let an integer $p\in\mathbb{Z}$, $\varepsilon\in\mathbb{Z}/2\mathbb{Z}$, we will define the sheaves of real analytic functions $\sha_{P}(\varepsilon,p)$, hyperfunctions $\shb_P(\varepsilon,p)$, microfunctions $\mathscr{C}_{P}(\varepsilon,p)$, on $P$ resp. $P^*$ twisted by some power of the tautological line bundle. 

For $X,Y$ either the manifold $\mathbb{P}$ or $\mathbb{P}^*$, for any two integers $p,q$, we note $\mathcal{O}_{X\times Y}(p,q)$ the line bundle on $X\times Y$ with homogenity $p$ in the $X$ variable and $q$ in the $Y$ variable. Setting $\Omega_{X\times Y/X}(p,q):=\Omega_{X\times Y/X} \tens_{\mathcal{O}_{X\times Y}}\mathcal{O}_{X\times Y}(p,q)$, $\she_X^{\mathbb{R}}(p,q):=\muhom(\mathbb{C}_{\Delta_{X}},\Omega_{X\times X/X}(p,q))[d_X]$ and we define accordingly $\she_X(p,q)$. Let us notice that $\she_X^{\mathbb{R}}(-p,p)$ is a sheaf of rings. 

Let $n$ be the dimension of $P$, (of course $n=d_{\mathbb{P}}$). For an integer $k$ and $\varepsilon\in\mathbb{Z}/2\mathbb{Z}$, we note $k^{*}:=-n-1-k$, $\varepsilon^{*}:=-n-1-\varepsilon\hspace{1ex}mod(2)$. We have:

\begin{theorem}\label{pre_main_theorem}
	\bnum
	\item Let $k$ be an integer such that $-n-1<k<0$ and let $s$ be a global non-degenerate section on $\dT{}^*_{\mathbb{S}}(\mathbb{P}\times\mathbb{P}^*)$ of $H^{1}(\muhom(\mathbb{C}_{\mathbb{S}},\Omega_{\mathbb{P}\times \mathbb{P}^*/\mathbb{P}^*}(-k,k^*)))$. For $P\in\she_\mathbb{P}(-k,k)$, there is a unique $Q\in\she_{\mathbb{P}^*}(-k^*,k^*)$ satisfying $P\cdot s=s\cdot Q$.
	The morphism induced by $s$
	\eqn
	\oim{\chi}\she_{\mathbb{P}}(-k,k)\to\she_{\mathbb{P}^*}(-k^*,k^*)
	\eneqn
	\vspace{-2em}
	\eqn
	P\mapsto Q
	\eneqn
	is a ring isomorphism.
	\item There exists such a non-degenerate section $s$.
	\enum
\end{theorem}

In fact, we will see that the non-degenerate section of Theorem \ref{pre_main_theorem} is provided by the Leray section.

Now, from classical adjunction formulas for $\she$-modules, we get a correspondance between solutions of systems of microdifferential equations on the projective space and solutions of systems of microdifferential equations on its dual. We will prove the following theorem, which was proved in \cite{DS96} for $\mathcal{D}$-modules,

\begin{theorem}
	Let $k$ be an integer such that $-n-1<k<0$. Let $\mathcal{N}$ be a coherent $\she_{\mathbb{P}}(-k,k)$-module and $F\in\Derb(\mathbb{P})$. Then, we have an isomorphism in $\Derb(\C_{\dT{}^*\mathbb{P}^*})$
	\eqn
	\hspace{-15em}\oim{\chi}\rhom{_{\she_{\mathbb{P}}(-k,k)}}(\mathcal{N},\muhom(F,\mathcal{O}_\mathbb{P}(k)))\simeq
	\eneqn
	\vspace{-2em}
	\eqn
	\hspace{15em}\rhom{_{\she_{\mathbb{P}^*}(-k^*,k^*)}}(\underline{\Phi}_{\mathbb{S}}^{\mu}(\mathcal{N}),\muhom((\Phi_{\mathbb{C}_{\mathbb{S}}[-1]}F,\mathcal{O}_{\mathbb{P}^*}(k^*)))
	\eneqn

\end{theorem}
where $\underline{\Phi}_{\mathbb{S}}^{\mu}$ it is the counterpart of $\underline{\Phi}_\mathbb{S}$ for $\she$-modules, and will be defined in Section \ref{sec:int_trans_emodules}.

Let us mention that, through a difficult result from \cite{KSIW06}, $\muhom(F,\mathcal{O}_\mathbb{P})$ is well-defined in the derived category of $\she$-modules.

\begin{corollary}\label{main_theorem}
	Let $k$ be an integer such that $-n-1<k<0$ and $\varepsilon\in\mathbb{Z}/2\mathbb{Z}$. The section $s$ of theorem \ref{pre_main_theorem} defines an isomorphism: 
	$$
	\chi_{*}\mathscr{C}_{P}(\varepsilon,k)\vert_{\dT{}^*_P\mathbb{P}}\simeq\mathscr{C}_{P^{*}}(\varepsilon^*,k^{*})\vert_{\dT{}^*_{P^*}\mathbb{P}^*}
	$$
	Moreover, this morphism is compatible with the respective action of $\oim{\chi}\she_{\mathbb{P}}(-k,k)$ and $\she_{\mathbb{P}^*}(-k^*,k^*)$.
\end{corollary}

\textbf{Acknowledgements} I would like to express my gratitude to Pierre Schapira for suggesting me this problem and for his enlightening insights to which this work owes much.

\section{Reminders on Algebraic Analysis and complements}\label{sec:reminder_algebraic_analysis}

In this section, we recall classical results of Algebraic Analysis, with the exception of section \ref{complements_muhom}.

\subsection{Notations for manifolds}\label{not:12345}

\bnum

\item Let $M_i$ ($i=1,2,3$) be manifolds. For short, we write
$M_{ij}\eqdot M_i\times M_j$ ($1\leq i,j\leq3$),
$M_{123}=M_1\times M_2\times M_3$,
$M_{1223}=M_1\times M_2 \times M_2\times M_3$, etc.

\item $\delta_{M_i}\cl M_i\to M_i\times M_i$ denote the diagonal embedding, and
$\Delta_{M_i}$ the diagonal set of $M_i\times M_i$.

\item We will often write for short $\cor_i$ instead of $\cor_{{M_i}}$ and $\cor_{\Delta_i}$
instead of $\cor_{\Delta_{M_i}}$
and similarly with $\omega_{M_i}$, etc.,
and with the index $i$ replaced with several indices $ij$,  etc.

\item We denote by $\pi_i$, $\pi_{ij}$, etc.\ the projection
$T^*M_{i}\to M_{i}$,
$T^*M_{ij}\to M_{ij}$, etc.

\item For a fiber bundle $E\to M$, we denote by $\dot{E}\to M$ the fiber bundle with the zero-section removed.

\item
We denote by $q_i$ the
projection $M_{ij}\to M_i$ or the projection $M_{123}\to M_i$ and by $q_{ij}$
the projection $M_{123}\to M_{ij}$. Similarly, we denote by $p_i$ the
projection $T^*M_{ij}\to T^*M_i$ or the projection $T^*M_{123}\to T^*M_i$ and
by $p_{ij}$ the projection $T^*M_{123}\to T^*M_{ij}$.

\item We also need to
introduce the maps $p_{j^a}$ or $p_{ij^a}$, the composition of $p_{j}$ or $p_{ij}$ and the antipodal
map $a$ on $T^*M_j$.  For example,
\eqn
&&p_{12^a}((x_1,x_2,x_3;\xi_1,\xi_2,\xi_3))=(x_1,x_2;\xi_1,-\xi_2).
\eneqn
\item
We let $\delta_2\cl M_{123} \to M_{1223}$ be the natural
diagonal embedding.
\enum

\medskip

\subsection{Sheaves}

We follow the notations of~\cite{KS90}. 

Let $X$ be a good topological space, i.e. separated, locally compact, countable at infinity, of finite global cohomological dimension and let $\cor$ be a commutative unital ring of finite global dimension. 

For a locally closed subset $Z$ of $X$, we denote by $\cor_Z$, the sheaf, constant on $Z$ with stalk $\cor$, and $0$ elsewhere. 

We denote by $\Derb(\cor_X)$ the bounded derived category of the category of sheaves of $\cor$-modules on $X$. If $\shr$ is a sheaf of rings, we denote by $\Derb(\shr)$ the bounded derived category of the category of left $\shr$-modules.

Let Y be a good topological space and $f$ a morphism $Y\to X$. We denote by $\roim{f},\opb{f},\reim{f},\epb{f},\rhom,\ltens$ the six Grothendieck operations. We denote by $\etens$ the exterior tensor product.

We denote by $\omega_X$ the dualizing complex on $X$, by $\omega_X^{\otimes-1}$ the sheaf-inverse of $\omega_X$ and by $\omega_{Y/X}$ the relative dualizing complex. 

 In the following, we assume that $X$ is a real manifold. Recall that $\omega_X\simeq\ori_X\,[\dim X]$ where $\ori_X$ is the orientation sheaf and $\dim X$ is the dimension of $X$. We denote by $\RD_X(\scbul)$,$\RD'_X(\scbul)$ the duality functor $\RD_X(\scbul)=\rhom(\scbul,\omega_X), \RD'_X(\scbul)=\rhom(\scbul,\cor_X)$, respectively.

For $F\in\Derb(\cor_X)$, we denote by $SS(F)$ its singular support, also called micro-support. For a a subset $Z\subset T^*X$, we denote by $\Derb(\cor_X;Z)$ the localization of the category $\Derb(\cor_X )$ by the full subcategory of objects whose micro-support is contained in $T^*X\setminus Z$.

For a closed submanifold $M$ of $X$, we denote by $\nu_M$, $\mu_M$, $\muhom$, the functor of specialization, microlocalization along $M$ and the functor of microcalization of $\rhom$ respectively.


 Let $M_i$ ($i=1,2,3$) be manifolds. We shall consider the operations of composition of kernels:
\eq\label{eq:conv}
&&\ba{l}
\conv[2]\;\cl\;\Derb(\cor_{M_{12}})\times\Derb(\cor_{M_{23}})\to\Derb(\cor_{M_{13}})\\
\hs{10ex}\ba{rcl}(K_1,K_2)\mapsto K_1\conv[2] K_2&\eqdot&
\reim{q_{13}}(\opb{q_{12}}K_1\ltens\opb{q_{23}}K_2)\\
&\simeq&\reim{q_{13}}\opb{\delta_2}(K_1\letens K_2)\ea
\ea
\eneq
\eq\hs{10ex}\label{eq:2_conv}
&&\ba{l}
\conv[23]\;\cl\;\Derb(\cor_{M_{12}})\times\Derb(\cor_{M_{23}})\times\Derb(\cor_{M_{34}})\to\Derb(\cor_{M_{14}})\\
\hs{10ex}\ba{rcl}(K_1,K_2,K_3)\mapsto K_1\conv[2] K_2 \conv[3] K_3 &\eqdot&
\reim{{q_{14}}}(\opb{{q_{12}}}K_1\ltens\opb{{q_{23}}}K_2\ltens\opb{{q_{34}}}K_3)\\
\ea
\ea
\eneq

Let us mention a variant of $\circ$:
\eqn
&&\ba{l}
\sconv[2]\;\cl\;\Derb(\cor_{M_{12}})\times\Derb(\cor_{M_{23}})
\to\Derb(\cor_{M_{13}})\\
\hs{10ex}(K_1,K_2)\mapsto K_1\sconv[2] K_2\eqdot
\roim{q_{13}}\bl\opb{q_{2}}\omega_{2}\tens\epb{\delta_2}(K_1\etens K_2)\br
\ea
\eneqn
There is a natural  morphism $K_1 \conv[2] K_2  \to K_1 \sconv[2] K_2$.

We refer the reader to \cite{KS90} for a detailed presentation of sheaves on manifolds.
%
%

\subsection{$\mathcal{O}$-modules and $\mathcal{D}$-modules}
We refer to \cite{K03} for the notations and the main results of this section.

Let $(X,\mathcal{O}_X)$ be a complex manifold. We denote by $d_X$ its complex dimension and by $\mathcal{D}_X$ the sheaf of rings of finite order holomorphic differential operators on $X$.

For an invertible $\mathcal{O}_X$-module $\mathcal{F}$, we denote by $\mathcal{F}^{\tens -1}:=\hom{_{\mathcal{O}_X}}(\mathcal{F},\mathcal{O}_X)$,
the inverse of $\mathcal{F}$. Denote by $\text{Mod}(\mathcal{D}_X)$ the abelian category of left $\mathcal{D}_X$-modules and $\text{Mod}(\mathcal{D}_X^{op})$ of right $\mathcal{D}_X$-modules. We denote by $\Omega_X$ the right $\mathcal{D}_X$-module of holomorphic $d_X$ forms. 
%

Let $\Derb(\mathcal{D}_X)$ be the bounded derived category of the category of left $\mathcal{D}_X$-modules, $\Derb_{\text{coh}}(\mathcal{D}_X)$ its full triangulated subcategory whose objects have coherent cohomology. 

%

Let $\Derb_{\text{good}}(\mathcal{D}_X)$ be the triangulated subcategory of $\Derb(\mathcal{D}_X)$, whose objects have all cohomologies consisting in good $\mathcal{D}_X$-modules (see \cite{K03} for a classical reference).


We refer in the following to \cite{K03}. Let $f:Y\to X$ be a morphism of complex manifolds. We denote by $\mathcal{D}_{Y\to X}$ and $\mathcal{D}_{X\from Y}$ the transfer bimodules.

For $\mathcal{M}\in\Derb(\mathcal{D}_X)$, $\mathcal{N}\in\Derb(\mathcal{D}_Y)$, we denote by $\opb{\underline{f}}\mathcal{M}$, $\oim{\underline{f}}\mathcal{N}$, the pull-back and the direct image of $\mathcal{D}$-modules respectively.
We refer to \cite{DS96} for functorial properties of inverse and direct image of $\mathcal{D}$-modules.

\subsection{$\she$-modules}\label{sec:reminder_e_mod}

We refer in the following to \cite{SKK73} (see also \cite{S85} for an exposition). For a complex manifold $X$, one denotes by $\she_X$ the sheaf of filtered ring of finite order holomorphic microdifferential operators on $T^*X$. We denote by $\Derb_{\text{coh}}(\she_X)$ the full triangulated subcategory of $\Derb(\she_X)$  whose objects have coherent cohomology. 

For $m\in\mathbb{Z}$, we denote by $\she_{X}(m)$ the abelian subgroup of $\she_X$ of microdifferential operators of order less or equal to $m$. For a section $P$ of $\she_X$, we denote by $\sigma(P)$ the principal symbol of $P$. 

Let $\pi_X$ denote the natural projection $T^*X\to X$. Let us recall that $\she_X$ is flat over $\opb{\pi}(\mathcal{D}_X)$. To a $\mathcal{D}_X$-module $\shm$, we associate an $\she_X$-module defined by 
\eqn
\she\shm:=\she_X\tens_{\opb{\pi_X}\mathcal{D}_X}\opb{\pi_X}\shm
\eneqn
 To a morphism of manifolds $f\cl Y\to X$, we associate the diagram of natural morphisms:
 \eq\label{diag:microlocal1}
 &&\xymatrix{
 	T^*Y\ar[dr]_-{\pi_Y}&\ar[l]_-{f_d}Y\times_XT^*X\ar[r]^-{f_\pi}\ar[d]^-\pi&T^*X\ar[d]^-{\pi_X}\\
 	&Y\ar[r]^-f&X
 }\eneq
 where $f_d$ is the transposed of the tangent map $Tf\cl TY\to Y\times_XTX$.

 For $\mathcal{M},\mathcal{N}$ objects of respectively $\Derb(\she_X)$ and $\Derb(\she_Y)$,  we denote by $\opb{\underline{f}}\mathcal{M}$ and $\oim{\underline{f}}\mathcal{N}$  the pull-back and the direct image of $\she$-modules respectively.

\subsection{Hyperfunctions and microfunctions}\label{sec:reminder_special_sheaves}

Let $M$ be a real analytic manifold and $X$ a complexification of $M$. We might be led to identify $T_M^*X$ with $i\cdot T^*M$. We denote by $\sha_M\eqdot\sho_X\vert_{M}$, $\shb_M\eqdot\rhom(\RD'_X\C_M,\sho_X)$, $\shc_M\eqdot\muhom(\RD'_X\C_M,\sho_X)$,
the sheaves of real analytic functions, hyperfunctions, microfunctions, respectively. Let us denote by \textit{sp}, the isomorphism
\eq\label{def:spectrum_isomorphism}
sp\cl \shb_M \isoto \roim{\pi_{M}}\shc_M
\eneq 


There is a natural action of the sheaf of microdifferential operators $\she_X$ on $\shc_M$.


If $Z$ is a closed complex submanifold of $X$ of codimension $d$, we note 
\eqn
\shb_{Z\vert X}:=H^d_{[Z]}(\mathcal{O}_X)
\eneqn
the algebraic cohomology of $\mathcal{O}_X$ with support in $Z$.

\subsection{Integral transforms for sheaves and $\mathcal{D}$-modules}\label{sec:d_modules_duality_intro}

\subsubsection{Integral transforms for sheaves}

Let $X$ and $Y$ be complex manifolds of respective dimension $d_X,d_Y$. Let $S$ be a closed submanifold $X\times Y$ of dimension $d_S$. We set $d_{S/X}:=d_S-d_X$. Consider the diagram of complex manifolds
\eq\label{eq:duality_diagram}
\xymatrix@C-0pc@R+0pc{
	&S \ar[dl]_-{f} \ar[dr]^-{g} &&&\widetilde{S} \ar[dl]_-{g} \ar[dr]^-{f} &\\  
	X & & Y,&  Y & & X
}
\eneq
where the second diagram is obtained by interchanging $X$ and $Y$.

Let $F\in\Derb(\mathbb{C}_X)$, $G\in\Derb(\mathbb{C}_Y)$, we define
\eqn
\xymatrix@C-0pc@R+0pc{
	\Phi_S(F):=\reim{g}\opb{f}F[d_{S/Y}],&\Phi_{\widetilde{S}}(G):=\reim{f}\opb{g}G[d_{S/X}]\\  
}
\eneqn
\eqn
\xymatrix@C-0pc@R+0pc{
	\Psi_S(F):=\roim{g}\epb{f}F[d_{X/S}],&\Psi_{\widetilde{S}}(G):=\roim{f}\epb{g}G[d_{Y/S}]\\  
}
\eneqn

For $K\in\Derb(\mathbb{C}_{X\times Y})$, and given the diagram $X\from[q_1] X\times Y \to[q_2] Y$, we define the integral transform of $F$ with kernel $K$
\eqn
\xymatrix@C-0pc@R+0pc{
	\Phi_K(F):=\reim{q_2}(K\tens\opb{q_1}F)\\  
}
\eneqn

%

\subsubsection{Integral transforms for $\mathcal{D}$-modules}

Let $X,Y$ be complex manifolds of equal dimension $n>0$, and $S$ a complex manifold. Consider again the situation (\ref{eq:duality_diagram}).

We suppose
\eq\label{hyp:duality_diagram}
\left\{
\begin{array}{ll}
	\hspace{1em}\text{$f,g$ are smooth and proper,}\\
	\hspace{1em}\text{$S$ is a complex submanifold of $X\times Y$ of codimension $c>0$}\\
\end{array}
\right.
\eneq

Let $\mathcal{M}\in\Derb(\mathcal{D}_X)$, $\mathcal{N}\in\Derb(\mathcal{D}_Y)$. Let us denote by $\widetilde{S}$ the image of $S$ by the map $r:X\times Y\to Y\times X, (x,y)\mapsto(y,x)$. One sets
\eqn
\xymatrix@C-0pc@R+0pc{
	\underline{\Phi}_S(\mathcal{M}):=\oim{\underline{g}}\opb{\underline{f}}\mathcal{M},&\underline{\Phi}_{\widetilde{S}}(\mathcal{N}):=\oim{\underline{f}}\opb{\underline{g}}\mathcal{N}
}
\eneqn

We refer to \cite[Prop. 2.6.]{DS94} for adjonction formulae related to these integral transforms.

Let us recall that we denote by $\Omega_X$ the sheaf of holomorphic  $n$-forms and let
\eqn
\shb_{S\vert X\times Y}^{(n,0)}:=\opb{q_1}\Omega_X\tens_{\opb{q_1}\mathcal{O}_X}\shb_{S\vert X\times Y}
\eneqn
This $(\mathcal{D}_{Y},\mathcal{D}_{X})$-bimodule allows the computation of $\underline{\Phi}_S$ because of the isomorphism, proven in \cite[Prop 2.12]{DS94}
\eqn
\mathcal{D}_{Y\from S}\ltens{_{\mathcal{D}_{S}}}\mathcal{D}_{S\to X}\isoto\shb_{S\vert X\times Y}^{(n,0)}
\eneqn

leading to
\eqn
\underline{\Phi}_S(\mathcal{M})\simeq\reim{q_2}(\shb_{S\vert X\times Y}^{(n,0)}\ltens{_{\opb{q_1}\mathcal{D}_{X}}}\opb{q_1}\mathcal{M})
\eneqn

\subsection{Microlocal integral transforms}

\subsubsection{Integral transforms for $\she$-modules}\label{sec:int_trans_emodules}

Let $X,Y$ be complex manifolds and $S$ is a closed submanifold of $X\times Y$. We consider again the diagram (\ref{eq:duality_diagram}) under the hypothesis (\ref{hyp:duality_diagram}).

We define the functor
\eqn
\begin{array}{rrrl}
	\Derb(\mathcal{E}_X)\to \Derb(\mathcal{E}_Y)\text{, } \underline{\Phi}_S^{\mu}(\mathcal{M}):=\oim{\underline{g}}\opb{\underline{f}}\mathcal{M}\\ 
\end{array}
\eneqn

We define the $\she_{X \times Y}$-module attached to $\shb_{S\vert X\times Y}$,
\eqn
\shc_{S\vert X\times Y}:=\she\shb_{S\vert X\times Y}
\eneqn

and we consider the $(\she_Y,\she_X)$-bimodule 
\eq\label{microlocal_associated_kernel}
\shc_{S\vert X\times Y}^{(n,0)}:=\opb{\pi}\opb{q_1}\Omega_X\ltens{_{\opb{\pi}\opb{q_1}\mathcal{O}_X}}\shc_{S\vert X\times Y}
\eneq

One can notice that
\eqn
\she_{Y\from S}\ltens{_{\she_{S}}}\she_{S\to X}\isoto\shc_{S\vert X\times Y}^{(n,0)}
\eneqn 

and hence, we have 
\eq\label{eq:microlocal_integral_transform}
\underline{\Phi}_S^{\mu}(\mathcal{M})\simeq\reim{p_2^a}(\shc_{S\vert X\times Y}^{(n,0)}\ltens{_{\opb{p_1}\she_{X}}}\opb{p_1}\mathcal{M})
\eneq

Let $\mathcal{M}\in\Derb_{\text{good}}(\mathcal{D}_X)$. The functors $\underline{\Phi}_S^{\mu}$ and $\underline{\Phi}_S$ are linked through the following isomorphism in $\Derb(\mathbb{C}_{\dT{}^*Y})$, (see \cite{SS94})
\eq\label{eq:e_d_modules_transform}
\xymatrix@C-0pc@R+0pc{
	\she(\underline{\Phi}_S(\mathcal{M}))\simeq\underline{\Phi}_S^{\mu}(\she\mathcal{M})\\ 
}
\eneq

\subsubsection{Microlocal integral transform of the structure sheaf}
Consider two open subsets $U$ and $V$ of $T^*X$ and $T^*Y$, respectively and $\Lambda$ a closed complex Lagrangian submanifold of $U\times V^a$
\eq\label{eq:1b}
&&\xymatrix{
	&U\times V^a\ar[ld]_-{p_1}\ar[rd]^-{p_{2^a}}&\\
	T^*X\supset U&  &V\subset T^*Y
}\eneq

As detailed in Section 11.4 of \cite{KS90}, let $K\in\Derb(\mathbb{C}_{X\times Y})$, $SS(K)$ its micro-support and let us suppose that $p_1\vert{\Lambda},p_2^a\vert{\Lambda}$ are isomorphisms, $K$ is cohomologically constructible simple with shift $0$ along $\Lambda$ and that $(\opb{p_1}(U)\cup\opb{p_2}(V))\cap SS(K)\subset\Lambda$.

Let $p=(p_X,p_Y^a)\in\Lambda$ and let us consider some section $s\in H^0(\muhom(K,\Omega_{X\times Y/Y}))_{p}$, where $\Omega_{X\times Y/Y}:=\mathcal{O}_{X\times Y}\tens_{\opb{q_1}\mathcal{O}_{X}}\opb{q_1}\Omega_X$. The section $s$ gives a morphism $K\to\Omega_{X\times Y/Y}$ in $\Derb(\mathbb{C}_{X\times Y};p)$. Then, there is a natural morphism
\eq
\begin{array}{ll}
	\Phi_{K[d_X]}(\mathcal{O}_X) & \to\mathcal{O}_Y
\end{array}%
\eneq

We recall the result:
\begin{theorem}[{{\cite[Th.11.4.9]{KS90}}}]
	There exists $s\in H^0(\muhom(K,\Omega_{X\times Y/Y}))_{p}$ such that the associated morphism $\Phi_{K[d_X]}(\mathcal{O}_X)\to\mathcal{O}_Y$  is an isomorphism in the category $\Derb(\mathbb{C}_Y;p_Y)$. Moreover, this morphism is compatible with the action of microdifferential operators on $\mathcal{O}_X$ in $\Derb(\mathbb{C}_X;p_X)$ and the action of microdifferential operators on $\mathcal{O}_Y$ in $\Derb(\mathbb{C}_Y;p_Y)$
\end{theorem}

Also, we will make use of the following theorem proven in~\cite[Th.~7.2.1]{KS90}:

\begin{theorem}[{{\cite[Th.~7.2.1]{KS90}}}]\label{main_theorem_contact_muhom}
	Let $K\in\Derb(X\times Y)$
	and assume that 
	\item (i) K is cohomologically constructible
	\item (ii) $(p_{1}^{-1}(U)\cup (p_{2}^{a})^{-1}(V))\cap SS(K)\subset\Lambda$
	\item (iii) the natural morphism $\mathbb{C}_{\Lambda}\longrightarrow \mu\mathpzc{hom}(K,K)|_\Lambda$ is an isomorphism.
	\item Then for any $F_{1},F_{2}\in\Derb(X;U)$, the natural morphism
	
	$$
	\chi_{*}\mu\mathpzc{hom}(F_{1},F_{2})\longrightarrow\mu\mathpzc{hom}(\Phi_K(F_{1}),\Phi_K(F_{2}))
	$$
	is an isomorphism in $\Derb(Y;V)$.
\end{theorem}

\subsection{Complements on the functor $\muhom$}\label{complements_muhom}

\subsubsection{Associativity for the composition of kernels}
The next result is well-known although no proof is written down in the literature, to our knowledge.

\begin{lemma}\label{eq:associativity_lemma}
	
	Let $M_{1},M_{2},M_{3}$ be real manifolds, and $K,L,M$ be objects respectively of $\Derb(\cor_{M_{12}})$, $\Derb(\cor_{M_{23}}),\Derb(\cor_{M_{34}})$, then the composition of kernels $\conv$ defined in \ref{eq:conv} is associative. We have the following isomorphism 
	
	\eq\label{diag:associativity_composition}
	\begin{array}{rcl}
	(K\conv[2] L)\conv[3] M & \isoto & K\conv[2] (L\conv[3] M)\\
	\end{array}	
	\eneq
	
	such that for any $N\in  \Derb(\cor_{M_{45}})$, the diagram below commutes:
	
	\eq\label{diag:tensor}
	&&\xymatrix{
	((K\conv[2] L)\conv[3] M) \conv[4] N\ar[r]\ar[d]&(K\conv[2] L)\conv[3] (M \conv[4] N)\ar[dd]\\
	(K\conv[2] (L\conv[3] M)) \conv[4] N\ar[d]&\\
	K\conv[2] ((L\conv[3] M )\conv[4] N)\ar[r]&K\conv[2] (L\conv[3] (M \conv[4] N)).
	}\eneq
	
\end{lemma}
\begin{proof}
	Consider the following diagram
	
	\eqn
	\xymatrix@C-1pc@R+1pc{
	&&&&&M_{1234}\ar@<0.0ex>[ld]_-{q_{124}^{3}}\ar@<0.15ex>[ld]\ar[dd]_-{q_{14}^{23}} \ar@<0.0ex>[rd]^-{q_{134}^{2}}\ar@<0.15ex>[rd]
	\ar@/_3pc/@<0.0ex>[lllldd]\ar@/_3pc/@<0.15ex>[lllldd]\ar@/^3pc/@<0.0ex>[rrrrdd]\ar@/^3pc/@<0.15ex>[rrrrdd]&&&&&\\
	&&&&M_{124}\ar@<0.0ex>[rd]\ar@<0.15ex>[rd]_-{q_{14}^{2}}\ar@/^1pc/@<0.0ex>[rrrrrdd]\ar@/^1pc/@<0.15ex>[rrrrrdd]^-{q_{24}^{1}}\ar@/_3pc/[lllldd]_-{q_{12}^{4}}&&
	M_{134}\ar@<0.0ex>[ld]^-{q_{14}^{3}}\ar@<0.15ex>[ld]\ar@/_1pc/@<0.0ex>[llllldd]_-{q_{13}^{4}}\ar@/_1pc/@<0.15ex>[llllldd]\ar@/^3pc/[rrrrdd]^-{q_{34}^{1}}&&&&\\
	&M_{123}\ar[ld]_-{q_{12}^{3}}\ar@<0.0ex>[d]^-{q_{13}^{2}}\ar@<0.11ex>[d]\ar[rrrrd]^-{q_{23}^{1}}&&&&M_{14}&&&&
	M_{234}\ar[lllld]_-{q_{23}^{4}}\ar@<0.0ex>[d]_-{q_{24}^{3}}\ar@<0.11ex>[d]\ar[rd]^-{q_{34}^{2}}&\\
	M_{12}&M_{13}&&&&M_{23}&&&&M_{24}&M_{34}
	}
	\eneqn
	
	where the thick squares are cartesian, and where for clarity we enforced the notation: the projection $M_{ijk}\to M_{ij}$ by $q_{ij}^{k}$ (independently of order of appearence of the indices), and the projection $M_{ijkl}\to M_{ij}$ by $q_{ij}^{kl}$. We now have:
	\eqn
	\begin{array}{rcl}
	\reim{{q_{14}^{3}}}(\opb{{q_{13}^{4}}}(\reim{{q_{13}^{2}}}(\opb{{q_{12}^{3}}}K\tens\opb{{q_{23}^{1}}}L))\tens\opb{{q_{34}^{1}}}M) & \simeq &  \reim{{q_{14}^{3}}}(\reim{{q_{134}^{2}}}(\opb{{q_{12}^{34}}}K\tens\opb{{q_{23}^{14}}}L)\tens\opb{{q_{34}^{1}}}M)\\
	& \simeq & \reim{{q_{14}^{23}}}(\opb{{q_{12}^{34}}}K\tens\opb{{q_{23}^{14}}}L\tens\opb{{q_{34}^{12}}}M)\\ 
	& \eqdot &
	K_1\conv[2] K_2 \conv[3] K_3
	\end{array}
	\eneqn
	
	The same way, we get the  isomorphism 
	\eqn
	\begin{array}{rcl}
	K_1\conv[2] K_2 \conv[3] K_3 & \simeq & \reim{{q_{14}^{2}}}(\opb{{q_{12}^{4}}}K \tens (\opb{{q_{24}^{1}}}(\reim{{q_{24}^{3}}}(\opb{{q_{23}^{4}}}L\tens\opb{{q_{34}^{2}}}M))))
	\end{array}
	\eneqn
	which proves the isomorphism (\ref{diag:associativity_composition}). And, it follows immediately that given $N\in  \Derb(\cor_{M_{45}})$, the diagram (\ref{diag:tensor}) commutes.

\end{proof}

%

\subsubsection{Associativity for the composition of $\muhom$}\label{sec:asso_composition_muhom}

We define the composition of kernels on cotangent bundles
(see~\cite[section 3.6, (3.6.2)]{KS90}).
\eq\label{eq:aconv}
&&\hs{-0ex}\ba{rcl}
\aconv[2]\;\cl\;\Derb(\cor_{T^*M_{12}})\times\Derb(\cor_{T^*M_{23}})
&\to&\Derb(\cor_{T^*M_{13}})\\
(K_1,K_2)&\mapsto&K_1\aconv[2] K_2\eqdot
\reim{p_{13}}(\opb{p_{12^a}} K_1\tens\opb{p_{23}} K_2)\\
&&\hs{8ex}\simeq\reim{p_{13^a}}(\opb{p_{12^a}} K_1
\tens\opb{p_{23^a}} K_2).
\ea
\eneq

There is a variant of the composition $\circ$, constructed in \cite{KS14}:
\eq\label{eq:star}
&&\ba{l}
\sconv[2]\;\cl\;\Derb(\cor_{M_{12}})\times\Derb(\cor_{M_{23}})
\to\Derb(\cor_{M_{13}})\\
\hs{10ex}(K_1,K_2)\mapsto K_1\sconv[2] K_2\eqdot
\roim{q_{13}}\bl\opb{q_{2}}\omega_{2}\tens\epb{\delta_2}(K_1\letens K_2)\br.
\ea
\eneq
There is a natural morphism for $K_1\in \Derb(\cor_{M_{12}})$ and $K_2\in \Derb(\cor_{M_{23}})$, $K_1 \conv[2] K_2  \to K_1 \sconv[2] K_2$.

Let us state a theorem proven in \cite[Prop.~4.4.11]{KS90} refined in \cite{KS14}. 

\begin{theorem}\label{th:main_associativity_theorem}
	Let $F_i,G_i,H_i$ respectively in $\Derb(\cor_{M_{12}}),\Derb(\cor_{M_{23}}),\Derb(\cor_{M_{34}})$, $i=1,2$. Let $U_i$ be an open subset of $T^*M_{ij}$ ($i=1,2$, $j=i+1$) and set $U_3=U_1\aconv[2]U_2$. There exists a canonical morphism in $\Derb(\cor_{T^*M_{13}})$, functorial in $F_1$ (resp. $F_2$):
	\eq\label{eq:micro_action_star_compo}
	&&\muhom(F_1,F_2)\vert_{U_1}\aconv[2]\muhom(G_1,G_2)\vert_{U_2}\to\muhom(F_1\sconv[2]G_1,F_2\conv[2]G_2)\vert_{U_3}.
	\eneq
	and hence
	\eq\label{eq:micro_action_compo}
	&&\muhom(F_1,F_2)\vert_{U_1}\aconv[2]\muhom(G_1,G_2)\vert_{U_2}\to\muhom(F_1\conv[2]G_1,F_2\conv[2]G_2)\vert_{U_3}.
	\eneq
\end{theorem}

We state the main theorem of this section.
\begin{theorem}\label{micro_associativity}
	Let $F_i,G_i,H_i$ respectively in $\Derb(\cor_{M_{12}}),\Derb(\cor_{M_{23}}),\Derb(\cor_{M_{34}})$, $i=1,2$ then we have:
	\banum
	\item
	\eqn
	\hspace{-15em}\left(\muhom(F_{1},F_{2})\aconv[2]\muhom(G_{1},G_{2})\right)\aconv[3]\muhom(H_{1},H_{2})\isoto
	\eneqn
	\vspace{-2em}
	\eqn
	\hspace{15em}\muhom(F_{1},F_{2})\aconv[2]\left(\muhom(G_{1},G_{2})\aconv[3]\muhom(H_{1},H_{2})\right)
	\eneqn

	\item The above isomorphism is compatible with the composition $\conv$ in the sense that the following diagram commutes
	\eqn
	\xymatrix@C-6pt{
	\hspace{-5em}(\muhom(F_{1},F_{2})\aconv[2]\muhom(G_{1},G_{2}))\aconv[3]\muhom(H_{1},H_{2})
	\ar[r]^-{\sim} \ar[d] &
	\muhom(F_{1},F_{2})\aconv[2](\muhom(G_{1},G_{2})\aconv[3]\muhom(H_{1},H_{2}))
	\ar[d]
	\\
	\muhom(F_{1}\conv[2] G_{1},F_{2}\conv[2] G_{2})\aconv[3]\muhom(H_{1},H_{2})
	\ar[d] &
	\muhom(F_{1},F_{2})\aconv[2]\muhom(G_{1}\conv[3] H_{1},G_{2}\conv[3] H_{2})
	\ar[d]
	\\
	\muhom((F_{1}\conv[2] G_{1})\conv[3] H_{1},(F_{2}\conv[2]G_{2})\conv[3] H_{2})
	\ar[r]^-{\sim} &
	\muhom(F_{1}\conv[2] (G_{1}\conv[3] H_{1}),F_{2}\conv[2] (G_{2}\conv[3] H_{2}))
	}
	\eneqn
	\eanum
	
\end{theorem}
\begin{proof}
	\banum
	\item
	This is a direct application of Lemma $\ref{eq:associativity_lemma}$ with $X,Y,Z$ taken to be respectively $T^*M_{12},T^*M_{13},T^*M_{34}$.
	
	\item We shall skip the proof, which is tedious but straightforward.
	\eanum
\end{proof}	
	
\section{Complex quantized contact transformations}

\subsection{Kernels on complex manifolds}\label{sec:kernels_complex_manifolds}

Consider two complex manifolds $X$ and $Y$ of respective dimension $d_X$ and $d_Y$.  We shall follow the notations of Section~\ref{not:12345}. 

For $K\in\Derb(\mathbb{C}_{X\times Y})$, we recall that we defined the functor $\Phi_K:\Derb(\mathbb{C}_Y)\rightarrow\Derb(\mathbb{C}_X)$, $\Phi_K(G)=Rq_{1 !}(K\tens q_{2}^{-1}(G)),\text{ for } G\in\Derb(\mathbb{C}_Y)$. With regards to the notation of Section~\ref{not:12345}, let us notice that $\Phi_K(G)$ is $K\circ G$. We refer also to Section \ref{The problem} for a definition of $\Omega_{X\times Y/X}$.

We recall the

\begin{lemma}\label{le:sects1}
	There is a natural  morphism
	\eqn
	&&\Omega_{X\times Y/X}\conv\sho_Y\,[d_Y]\to\sho_X.
	\eneqn
\end{lemma}
\begin{proof}
	We have
	\eqn
	\Omega_{X\times Y/X}\conv\sho_Y\,[d_Y]&=&
	\reim{q_1}(\Omega_{X\times Y/X}\tens\opb{q_1}\sho_Y[d_Y])\\
	&\to&\reim{q_1}(\Omega_{X\times Y/X}[d_Y])\to[\int]\sho_X,
	\eneqn
	where the last arrow is
	the integration morphism on complex manifolds.
\end{proof}

The following Lemma will be useful for the proof of Lemma \ref{lem:conv_section}. 
Let us first denote by $M_i$ ($i=1,2,3,4$) four complex manifolds, $L_i\in\Derb(\C_{M_{i,i+1}})$, $1\leq i\leq3$.
We set for short
\eqn
&&d_i=\dim_{\C} M_i, d_{ij}=d_{i}+d_{j}, \Omega_{ij/i}=\Omega_{M_{ij}/M_i}=\Omega^{(0,d_j)}_{M_{ij}}.
\eneqn
Set for $1\leq i\leq3, $
\eqn
&&K_{i}=\muhom(L_{i},\Omega_{i,j/i}[d_j]),\quad j=i+1\\
&&L_{ij}=L_i\circ L_j\quad  j=i+1,\quad L_{123}=L_1\circ L_2\circ L_3,\\
&&\tw K_{ij}=\muhom(L_{ij},\Omega_{i,j/i}\,[d_j]\circ \Omega_{j,k/j}\,[d_k]) \quad j=i+1, k=j+1\\
&&\tw K_{123}=\muhom(L_{123}, \Omega_{12/1}[d_2]\circ\Omega_{23/2}[d_3]\circ\Omega_{34/3}[d_4])\\
&&K_{ij}=\muhom(L_{ij},\Omega_{i,k/i}[d_{k}])\quad j=i+1, k=j+1,\\
&&K_{123}=\muhom(L_{123}, \Omega_{14/1}[d_4]).
\eneqn

We recall that we have t\\
he sequence of natural morphisms:
\eq\label{compo_forms_integration}
\Omega_{i,j/i}\conv \Omega_{j,k/j}&=&
\reim{q_{i,k}}(\opb{q_{i,j}}\Omega_{i,j/i}\tens\opb{q_{j,k}}\Omega_{j,k/j}) \nonumber\\
&\to&\reim{q_{i,k}}(\Omega_{i,j,k/i}) \nonumber\\
&\to&\Omega_{i,k/i}[-d_{j}]
\eneq

\begin{lemma}\label{lem:conv_int_commutation}
	The following diagram commutes:
	\eqn
	&&\xymatrix{
	&K_1\circ K_2\circ K_3\ar[ld]\ar[rd]\ar@{}[d]|-A&\\
	\tw K_{12}\circ K_3\ar[r]\ar[d]&\tw K_{123}\ar[d]&K_1\circ\tw K_{23}\ar[l]\ar[d]\\
	K_{12}\circ K_3\ar[r]\ar@{}[ru]|-B&K_{123}&K_1\circ K_{23}\ar[l]\ar@{}[lu]|-C
	}\eneqn
\end{lemma}
\begin{proof}
	Diagram labelled A commutes by the associativity of the functor $\muhom$ (see Theorem 2.7.3). Let us prove that Diagram B and C  commute. Of course, it is enough to consider Diagram B. To make the notations easier, we assume that $M_1=M_4=\rmpt$. We are reduced to prove the commutativity of the diagram:
	\eqn
	\xymatrix{
	\muhom(L_{2},\Omega_{2}\circ \Omega_{2,3/2}[d_{23}])\circ \muhom(L_{3},\sho_3)\ar[r]\ar[d]^-{\int_2}
	&\muhom(L_{23},\Omega_{2}\circ \Omega_{2,3/2}[d_{23}]\circ\sho_3)\ar[d]^-{\int_2}\\
	\muhom(L_{2},\Omega_{3}[d_{3}])\circ \muhom(L_{3},\sho_3)\ar[r]&\muhom(L_{23},\Omega_{3}[d_{3}]\circ\sho_3)
	}
	\eneqn
	For $F,F'\in\Derb(\cor_{12})$, $G,G'\in\Derb(\cor_{23})$, we saw in Theorem \ref{micro_associativity} (b) that the morphism $\muhom(F,F')\conv\muhom(G,G')\to\muhom(F\conv G,F'\conv G')$ is functorial in $F,F',G,G'$. This fact applied to the morphism
	\eqn
	\Omega_{2}\circ \Omega_{2,3/2}[d_{23}] \to \Omega_{3}[d_{3}]
	\eneqn
	gives that the above diagram commutes and so diagram $B$ commutes.
\end{proof}

Let $Z$ be a complex manifold and let $\Lambda\subset T^*(X\times Y)$ and $\Lambda'\subset T^*(Y\times Z)$ be two conic Lagrangian smooth locally closed complex submanifolds.

Let $L$, $L'$, be perverse sheaves on $X\times Y$, $Y\times Z$, with microsupport $SS(L)\subset\Lambda$, $SS(L')\subset\Lambda'$ respectively. We set
\eqn
&&L''\eqdot L[d_Y] \circ L'
\eneqn

Assume that
\eq\label{set_transversality_condition}
p_2^{a}\vert_{\Lambda}:\Lambda\to T^*{Y} \text{ and } p_2\vert_{\Lambda'}:\Lambda'\to T^*{Y} \text{ are transversal}
\eneq
 and that
 \eq\label{set_composability_condition}\hspace{3em}
 \text{the map } \Lambda\times_{T^*Y}\Lambda'\to\Lambda\conv\Lambda'
 \text{ is an isomorphism.}
 \eneq 

Let us set
\eqn
&&\shl\eqdot\muhom(L,\Omega_{X\times Y/X})
\eneqn
Note that $\shl\in\Derb(T^*(X\times Y))$ is concentrated in degree $0$. Indeed, it is proven in ~\cite[Th.~10.3.12]{KS90} that perverse sheaves are the ones which are pure with shift zero at any point of the non singular locus of their microsupport. On the other hand, Theorem 9.5.2 of ~\cite{KS85} together with Definition 9.5.1 of \cite{KS85} show that the latter verify the property that, when being applied $\muhom(\bullet,\Omega_{X\times Y/X})$, they are concentrated in degree $0$. Moreover, $\shl$ is a $(\she_X,\she_Y)$-bimodule. Indeed, such actions come from morphism (\ref{eq:micro_action_compo}) and the integration morphism (\ref{compo_forms_integration}). We define similarly $\shl'$ and $\shl''$.

Now consider two open subsets $U$,$V$ and $W$ of 
$\dTX$, $\dTY$, $\sdot{T\mspace{2mu}}{}^*{Z}$, respectively.

Let $K_{U\times V^a}$ be the constant sheaf on $(U\times V^a)\cap \Lambda$ with stalk $H^0\rsect(U\times V^a;\shl)$, extended by $0$ elsewhere.

$K'_{V\times W^a}$ is the constant sheaf on $(V\times W^a)\cap \Lambda'$ with stalk $H^0\rsect(V\times W^a;\shl')$, extended by $0$ elsewhere.

$K''_{U\times W^a}$ is the constant sheaf on $(U\times W^a)\cap \Lambda\circ\Lambda'$ with stalk $H^0\rsect(U\times W^a;\shl'')$, extended by $0$ elsewhere.

Let $s,s'$ be sections of $\sect(U\times V^a;\shl)$ and $\sect(V\times W^a;\shl')$ respectively. We define the product $s\cdot s'$ to be the section of $\sect(U\times W^a;\shl'')$, image of $1$ by the following sequence of morphisms
\eqn
\begin{array}{rcl}
	\mathbb{C}_{\Lambda\circ\Lambda'}&\isofrom&\mathbb{C}_{\Lambda}\circ\mathbb{C}_{\Lambda'}\\
	&:=&
	\reim{p_{13}}(\opb{p_{12^a}}\mathbb{C}_{\Lambda}\tens \opb{p_{23}}\mathbb{C}_{\Lambda'})\\
	&\to&\reim{p_{13}}(\opb{p_{12^a}}K_{U\times V^a}\tens \opb{p_{23}}K_{V\times W^a})\\
	&\to&
	\reim{p_{13}}(\opb{p_{12^a}}\muhom(L,\Omega_{X\times Y/X})\tens \opb{p_{23}}\muhom(L',\Omega_{Y\times Z/Y}))\\
	&:=&\muhom(L,\Omega_{X\times Y/X})\circ \muhom(L',\Omega_{Y\times Z/Y})
	\to\shl''
\end{array}
\eneqn
where the first isomorphism comes from the assumption \ref{set_composability_condition}.

\begin{lemma}\label{lem:conv_section}
Assume that conditions \ref{set_transversality_condition} and \ref{set_composability_condition} are satisfied. Let $s,s'$ be sections of $\sect(U\times V^a;\shl)$ and $\sect(V\times W^a;\shl')$ respectively, and let $G\in\Derb(\mathbb{C}_Y)$, $H\in\Derb(\mathbb{C}_Z)$. Then,
\bnum
\item $s$ defines a morphism
\eqn
&&\alpha_G(s)\cl\mathbb{C}_{\Lambda}\circ\muhom(G,\sho_Y)\vert_V\to\muhom(L[d_Y]\conv G,\sho_X)\vert_U
\eneqn
\item Considering the morphism 
\eqn
&\alpha_H(s\cdot s')\cl\mathbb{C}_{\Lambda\circ\Lambda'}\circ\muhom(H,\sho_Z)\vert_W\to\muhom(L[d_Y]\circ L'[d_Z]\conv H,\sho_X)\vert_U
\eneqn
we have the isomorphism
\eqn
	\alpha_H(s\cdot s') \simeq \alpha_{L'[d_Z]\circ H}(s)\circ \Phi_{\mathbb{C}_{\Lambda}}(\alpha_H(s'))
\eneqn
\enum
\end{lemma}

\begin{proof}
	\bnum
	\item Given $s$ and two objects $G_1,G_2\in\Derb(\mathbb{C}_Y)$, we have a morphism 
	\eqn
	\mathbb{C}_{\Lambda}\circ\muhom(G_1,G_2)\vert_V&\to&\muhom(L\conv G_1,\Omega_{X\times Y/X}\conv G_2)\vert_U
	\eneqn
	corrresponding to the composition of morphisms:
	\eq\label{decomposition_morphism_s}
	\begin{array}{rcl}
	\reim{p_{1}}(\mathbb{C}_{\Lambda}\tens\opb{p_{2^a}}\muhom(G_1,G_2)\vert_V)&\to &\reim{p_{1}}(K_{U\times V^a}\tens\opb{p_{2^a}}\muhom(G_1,G_2)\vert_V)\\
	&\to&\reim{p_{1}}(\muhom(L,\Omega_{X\times Y/X})\tens\opb{p_{2^a}}\muhom(G_1,G_2)\vert_V)\\
	&\to&\muhom(L\conv G_1,\Omega_{X\times Y/X}\conv G_2)\vert_U
	\end{array}
	\eneq 
	
	where the second morphism comes from the natural morphism $K_{U\times V^a}\to\muhom(L,\Omega_{X\times Y/X})$. We conclude by choosing, $G_1=G$, $G_2=\mathcal{O}_Y$ and by using Lemma~\ref{le:sects1}:
	\eqn
	&&\xymatrix@R=0ex@C=0ex{
	\muhom(L\conv G_1,\Omega_{X\times Y/X}\conv \mathcal{O}_Y)&\to&\muhom(L\conv G_1, \mathcal{O}_X[-d_Y])&\isoto&\muhom(L[d_Y]\conv G_1, \mathcal{O}_X)
	}
	\eneqn

	\item  Let $H\in\Derb(\mathbb{C}_Z)$. We denote by $\shh\eqdot\muhom(H,\sho_Z)$. It suffices to prove that the following diagram commutes:
	\clearpage
	\eqn
	\xymatrix@C-1pc@R+1pc{\hspace{-3em}
	(\mathbb{C}_{\Lambda}\circ\mathbb{C}_{\Lambda'})\circ\shh \ar[r]^-{\simeq}\ar[d]_-{\simeq} \ar@{}[rdddddd]|-*+[o][F-]{A} \ar@{.>}[rdddd]^(.4){ \Phi_{\mathbb{C}_{\Lambda}}(\alpha_H(s'))} & \mathbb{C}_{\Lambda}\circ(\mathbb{C}_{\Lambda'}\circ\shh)\ar[d]\\
	\mathbb{C}_{\Lambda}\circ\mathbb{C}_{\Lambda'}\circ\shh \ar@{.>}[rddddddd]^(.4){\alpha(s\cdot s')} \ar[d]  & \mathbb{C}_{\Lambda}\circ K'\circ \shh \ar[d] \\
	K\circ K' \circ \shh \ar[ddd] & \mathbb{C}_{\Lambda}\circ \muhom(L',\Omega_{Y\times Z/Y}) \circ \shh \ar[d] \\
	& \mathbb{C}_{\Lambda}\circ\muhom(L'\circ H,\Omega_{Y\times Z/Y}\circ \sho_Z) \ar[d]^-{\int_Z} \\
	&  \mathbb{C}_{\Lambda}\circ\muhom(L'[d_Z]\circ H,\sho_Y) \ar[d]\ar@/_-3pc/@{.>}[dddd]^(.35){\alpha_{L'[d_Z]\circ H}(s)}  \\
	K\circ \muhom(L',\Omega_{Y\times Z/Y}) \circ \shh \ar[r]^(.6){\int_Z}\ar[d] \ar@{}[rd]|-*+[o][F-]{B}&K\circ \muhom(L'[d_Z]\circ H, \sho_Y) \ar[d] \\
	\muhom(L,\Omega_{X\times Y/X})\circ \muhom(L',\Omega_{Y\times Z/Y})\circ \shh \ar[r]^-{\int_Z} \ar[dd] \ar@{}[rdd]|-*+[o][F-]{C}  & \muhom(L,\Omega_{X\times Y/X}) \circ \muhom( L'[d_Z]\circ H, \sho_Y) \ar[d] \\
	& \muhom(L\circ L'[d_Z]\circ H,\Omega_{X\times Y/X} \circ \sho_Y) \ar[d]^-{\int_Y} \\
	\muhom(L\circ L' \circ H,\Omega_{X\times Y/X}\circ \Omega_{Y\times Z/Y}\circ\sho_Z) \ar[r]^-{\int_{Y,Z}}& \muhom(L[d_Y]\circ L'[d_Z]\circ H,\sho_X)
	}
	\eneqn

where we omitted the subscript $U\times V^a$ and $V\times W^a$, $H$, $L'[d_Z]\circ H$ for $K_{U\times V^a}$, $K'_{V\times W^a}$, $\alpha_H, \alpha_{L'[d_Z]\circ H}$, respectively. 

We know from Theorem \ref{th:main_associativity_theorem} that the operation $\conv$ is functorial, so that diagram $A$ and $B$ commute. For instance, diagram $A$ decomposes this way: 
\eqn
\xymatrix@C-1pc@R+1pc{
	\mathbb{C}_{\Lambda}\circ\mathbb{C}_{\Lambda'}\circ\shh \ar[rr]^-{\simeq}\ar[d]  && \mathbb{C}_{\Lambda}\circ(\mathbb{C}_{\Lambda'}\circ\shh)\ar[d]\\
	\mathbb{C}_{\Lambda}\circ K' \circ \shh \ar[r] \ar[dd] & \mathbb{C}_{\Lambda}\circ \muhom(L',\Omega_{Y\times Z/Y}) \circ \shh \ar[r]^-{\simeq} \ar[dd] & \mathbb{C}_{\Lambda}\circ (\muhom(L',\Omega_{Y\times Z/Y}) \circ \shh) \ar[d]^-{\int_Y}\\
	&& \mathbb{C}_{\Lambda}\circ \muhom(L'[d_Z]\circ H,\sho_Y) \ar[d]\\
	K\circ K' \circ \shh \ar[r] & K\circ \muhom(L',\Omega_{Y\times Z/Y}) \circ \shh \ar[r]^-{\int_Y} & K \circ \muhom(L'[d_Z]\circ H,\sho_Y)
}
\eneqn

Besides, diagram $C$ commutes by Lemma \ref{lem:conv_int_commutation}. 

	Finally, the bottom diagonal punctured line correponds to $\alpha(s\cdot s')$, since the following diagram commutes 
	\eqn
	\xymatrix@C-1pc@R+1pc{
	\mathbb{C}_{\Lambda\circ\Lambda'} \ar[r]^-{\simeq}\ar@{.>}[rd]^(0.5){\alpha(s\cdot s')} \ar[d]& \mathbb{C}_{\Lambda}\circ\mathbb{C}_{\Lambda'} \ar[d]\\ 
	K''_{U\times W^a}\ar[r] & \muhom(L\circ L'[d_Y],\Omega_{X\times Z/X})
	}
	\eneqn
	
	\enum
\end{proof}
\begin{remark}
	In the following, unless necessary, we will omit the subsript for $\alpha$.
\end{remark}
\begin{theorem}\label{th:KS14b}
	Let $s\in \sect(U\times V^a;\shl)$, $G\in\Derb(\mathbb{C}_Y)$. Then, 
	
	(i) $s$ defines a morphism
	\eq\label{eq:mors2}
	&&\alpha(s)\cl\mathbb{C}_{\Lambda}\circ\muhom(G,\sho_Y)\vert_V\to\muhom(L[d_Y]\conv G,\sho_X)\vert_U.
	\eneq
	(ii) Moreover, if $P\in\sect(U;\she_X)$  and $Q\in\sect(V;\she_Y)$  satisfy $P\cdot s=s\cdot Q$, then the diagram below commutes
	\eq\label{diag:Ps=sQ}
	&&\xymatrix{
	\mathbb{C}_{\Lambda}\circ\muhom(G,\sho_Y)\ar[rr]^-{\alpha(s)}\ar[d]_-{\Phi_{\mathbb{C}_{\Lambda}}(\alpha(Q))}&&\muhom(L[d_Y]\conv G,\sho_X)\ar[d]^-{\alpha(P)}\\
	\mathbb{C}_{\Lambda}\circ\muhom(G,\sho_Y)\ar[rr]_-{\alpha(s)}&&\muhom(L[d_Y]\conv G,\sho_X).
	}\eneq
\end{theorem}

\begin{proof}

	(i) is already proven in Lemma \ref{lem:conv_section}.

	(ii) With regards to the notation of Lemma \ref{lem:conv_section}, we consider the triplet of manifolds $X,X,Y$, $\Lambda=\mathbb{C}_{\Delta_X}$, $\shl:= \muhom(\mathbb{C}_{\Delta_X}[-n],\Omega_{X\times X/X})$. Then, the assumption \ref{set_transversality_condition} is satisfied and noticing that $\Phi_{\mathbb{C}_{\Delta_X}}\simeq Id_X$, we conclude by Lemma \ref{lem:conv_section} that
	\eqn
	\alpha(P) \circ \alpha(s)\simeq \alpha(P\cdot s) \simeq \alpha(s\cdot Q) \simeq \alpha(s)\circ \Phi_{\mathbb{C}_{\Lambda}}(\alpha(Q))
	\eneqn
	
\end{proof}

\subsection{Main theorem}\label{sec:qct_main_theorem}
In this section, we will apply Theorem \ref{th:KS14b} when we are given a homogeneous symplectic isomorphism. Let us recall some useful results.

For $\shm$ a left coherent $\she_X$-module generated by a section $u\in\shm$, we denote by $\shi_{\shm}$ the annihilator left ideal of $\she_X$ given by: 
\eqn
\shi_{\shm}:=\{P\in\she_X;Pu=0\}
\eneqn 
and by $\overline{\shi}_{\shm}$ the symbol ideal associated to $\shi_{\shm}$: 
\eqn
\overline{\shi}_{\shm}:=\{\sigma(P);P\in\shi_{\shm}\}
\eneqn

\begin{definition}[{\cite{K03}}]
	Let $\shm$ be a coherent $\she_X$-module generated by an element $u\in\shm$. We say that $(\shm,u)$ is a simple $\she_X$-module if $\overline{\shi}_{\shm}$ is reduced and $\overline{\shi}_{\shm}=\{\phi\in\sho_{T^*X};\phi\vert_{supp(\shm)}=0\}$.
\end{definition}

Consider two complex manifolds $X$ and $Y$, open subsets $U$ and $V$ of $\dT{}^* X$ and $\dT{}^* Y$, respectively, and denote by $p_1$ and $p_2$ the projections $U \xleftarrow{p_1} U\times V^a \xrightarrow{p_2} V$. Let $\Lambda$ be a smooth closed submanifold Lagrangian of $U\times V^a$. We will make use of the following result from~\cite[Th.~4.3.1]{SKK73},~\cite[Prop.~8.5]{K03}:

\begin{theorem}[{\cite{SKK73}},{\cite{K03}}]\label{th:quantized_iso}
	Let $(\shm,u)$ be a simple $\she_{X\times Y}$-module defined on $U\times V^a$ such that $\supp{\shm}=\Lambda$. Assume $\Lambda\to U$ is a diffeomorphism. Then, there is an isomorphism of $\she_{X}$-modules: 
	\eqn
	\begin{array}{rrl}
	\she_X\vert_{U} & \isoto & \oim{(p_{1}\vert_{U\times V^a})}\shm\\
	P & \mapsto & P\cdot u
	\end{array}%
	\eneqn
\end{theorem}

Assume that the projections $p_1\vert_\Lambda$ and $p_2^a\vert_\Lambda$ induce isomorphisms. We denote by $\chi$ the homogeneous symplectic isomorphism $\chi:=p_2\vert_\Lambda\circ\opb{p_1\vert_\Lambda}$, 

\eq\label{eq:intro_contact_iso_diagram}\label{th:quantized_iso}
&&\xymatrix{
	&\Lambda\subset U\times V^a\ar[ld]_\sim^-{p_1\vert_\Lambda}\ar[rd]^\sim_-{p_2^a\vert_\Lambda}&\\
	\dTX\supset U\ar[rr]^-\sim_-\chi&  &V\subset \dTY
}
\eneq
\begin{corollary}\label{coro:quantized_iso_microdiff}
	Let $(\shm,u)$ be a simple $\she_{X\times Y}$-module defined on $U\times V^a$. Assume $\supp{\shm}=\Lambda$. Then, in the situation of $(\ref{eq:intro_contact_iso_diagram})$, we have an anti-isomorphism of algebras 
	\eqn
	\oim{\chi}\she_X\vert_U\simeq\she_Y\vert_V
	\eneqn
\end{corollary}

Consider two complex manifolds $X$ and $Y$ of the same dimension $n$, open subsets $U$ and $V$ of $\dT{}^*X$ and $\dT{}^*Y$, respectively, $\Lambda$ a smooth closed Lagrangian submanifold of $U\times V^a$ and assume that the projections $p_1\vert_\Lambda$ and $p_2^a\vert_\Lambda$ induce isomorphisms, hence a homogeneous symplectic isomorphism $\chi\cl U\isoto V$:

\eq\label{eq:contact_iso_diagram}
&&\xymatrix{
	&\Lambda\subset U\times V^a\ar[ld]_\sim^-{p_1}\ar[rd]^\sim_-{p_2^a}&\\
	\dTX\supset U\ar[rr]^-\sim_-\chi&  &V\subset \dTY
}
\eneq
We consider a perverse sheaf $L$ on $X\times Y$ satisfying
\eq\label{eq:condition_microsupport}
&&(\opb{p_1}(U)\cup\opb{{p_2^a}}(V))\cap\SSi(L)=\Lambda.
\eneq
and a section $s$ in $\sect(U\times V^a;\muhom(L,\Omega_{X\times Y/X}))$.

	Let $G\in\Derb(\mathbb{C}_Y)$. From Theorem \ref{th:KS14b} (i), the left composition by $s$ defines the morphism $\alpha(s)$ in $\Derb(\C_U)$:
	\eq
	&&\mathbb{C}_{\Lambda}\circ\muhom(G,\sho_Y)\vert_V\xrightarrow{\alpha(s)}\muhom(L[n]\conv G,\sho_X)\vert_U
	\eneq
	
	The condition (\ref{eq:condition_microsupport}) implies that $\supp(\muhom(L,\Omega_{X\times Y/X})\vert_{\opb{{p_2^{a}}}(V)})\subset\Lambda$. Since, $p_1$ is an isomorphism from $\Lambda$ to $U$ and that $\chi\circ p_1\vert_\Lambda=p_2^a\vert_\Lambda$, we get a morphism in $\Derb(\C_U)$
	\eq\label{eq:contact_transform_1}
	&&\opb{\chi}\muhom(G,\sho_Y)\vert_V\xrightarrow{\alpha(s)}\muhom(\Phi_{L[n]}(G),\sho_X)\vert_U
	\eneq
	
	\begin{theorem}\label{th:qct_main_theorem}
	Assume that the section $s$ is non-degenerate on $\Lambda$. Then, for $G\in\Derb(\C_Y)$, we have the following isomorphism in $\Derb(\C_U)$
	\eq\label{eq:qct_main_theorem}
	&&\opb{\chi}\muhom(G,\sho_Y)\vert_V\isoto\muhom(\Phi_{L[n]}(G),\sho_X)\vert_U
	\eneq
	Moreover, this isomorphism is compatible with the action of $\she_Y$ and $\she_X$ on the left and right side of (\ref{eq:qct_main_theorem}) respectively. 
	\end{theorem}
	\begin{proof}
	Let us first prove the following lemma, whose proof is available at the level of germs  in \cite[Th. 11.4.9]{KS90}.

	Let us prove that the morphism \eqref{eq:contact_transform_1} is an isomorphism. Let $L^*$ be the perverse sheaf $\opb{r}\rhom(L,\omega_{X\times Y/Y})$ where $r$ is the map $X\times Y\to Y\times X,(x,y)\mapsto(y,x)$. Let $s'$ be a section of $\muhom(L^*,\Omega_{Y\times X/Y})$, non-degenerate on $r(\Lambda)$, then we apply the same precedent construction to get a natural morphism
	\eqn
	\oim{\chi}\muhom(\Phi_{L[n]}(G),\sho_X)\vert_U\to\muhom(\Phi_{L^{*}[n]}\circ\Phi_{L[n]}G,\sho_Y)\vert_V\simeq &&\muhom(\Phi_{L^{*}\circ L[n]}G,\sho_Y)\vert_V
	\eneqn

	We know from~\cite[Th.~7.2.1]{KS90} that $\C_{\Delta_{X}}\simeq L^{*}\circ L$, so that we get a morphism in $\Derb(\C_V)$
	\eq\label{eq:contact_transform_2}
	&&\oim{\chi}\muhom(\Phi_{L[n]}(G),\sho_X)\vert_U\xrightarrow{\alpha(s')}\muhom(G,\sho_Y)\vert_V
	\eneq

	We must prove that (\ref{eq:contact_transform_1}) and (\ref{eq:contact_transform_2}) are inverse to each other. By Lemma \ref{lem:conv_section}{(ii)}, we get that the composition of these two morphisms is $\alpha(s'\cdot s)$, with $s'\cdot s\in\she_X$.

	For any left $\she_X$-module $\shm$, corresponds a right $\she_X$-module $\Omega_X\tens_{\sho_X}\shm$. Fixing a non-degenerate form $t_X$ of $\Omega_X\vert_U$ (resp.  $t_Y$ of $\Omega_Y\vert_{V}$), we apply now Theorem~\ref{th:quantized_iso}: $s$ and $s'$ are non-degenerate sections so that $(\she_{X\times Y}\vert_{U\times V^a},t_X\tens s)$ and $(\she_{Y\times X}\vert_{V^a\times U},s'\tens t_Y)$ are simple and so isomorphic to $(\opb{p_1}\she_X\vert_U,1)\simeq(\opb{{p_2^a}}\she_Y\vert_{V},1)$. $\Omega_X$ resp. $\Omega_Y$ being invertible $\sho_X$-module resp. $\sho_Y$-module, we get as well for the left-right $(\she_X\vert_U,\she_Y\vert_{V})$ bi-module, resp. left-right $(\she_Y\vert_{V},\she_X\vert_U)$ bi-module generated by $s$ resp. $s'$, that they are both isomorphic to $\opb{p_1}\she_X\vert_U\simeq\opb{{p_2^a}}\she_Y\vert_{V}$. 

	Then, following the proof of \cite[Th. 11.4.9]{KS90}, $s $ and $s'$, define ring isomorphisms associating to each $P\in\she_X(U)$,  $P'\in\she_X(U)$, some $Q\in\she_Y(V)$, $Q'\in\she_Y(V)$, such that $P\cdot s=s\cdot Q$, $s'\cdot P'=Q'\cdot s'$, respectively. Hence, we get that $\alpha(s')\circ \alpha(s)$ is an automorphism $\muhom(G,\sho_Y)\vert_V$, defined by the left action of $s'\cdot s\in\she_X$. Hence, we can choose $s'$ so that $\alpha(s')\circ \alpha(s)$ is the identity.

	We are now in a position to prove Theorem \ref{th:qct_main_theorem}: we constructed in the proof of the lemma, for each  $P\in\opb{p_{1}}\she_X\vert_U$, some $Q\in\opb{p_{2^a}}\she_Y\vert_{V}$ such that $P\cdot s=s\cdot Q$ and we can apply Theorem~\ref{th:KS14b} to conclude.
\end{proof}

\section{Radon transform for sheaves}\label{part:radon_transform}

We are going to apply the results of the last chapter to the case of projective duality. Recall projective duality for $\mathcal{D}$-modules were performed by D'Agnolo-Schapira \cite{DS96}. We will extend their results in a micorlocal setting.

\subsection{Notations}\label{sec:notations_projective}

In the following, we will quantize the contact transform associated with the Lagrangian submanifold $\dT{}_\mathbb{S}^*(\mathbb{P}\times\mathbb{P}^*)$, where $\mathbb{S}$ is the hypersurface of $\mathbb{P}\times \mathbb{P}^{*}$ defined by the incidence relation $\langle \xi,x\rangle=0,(x,\xi)\in\mathbb{P}\times \mathbb{P}^{*}$.

We denote by $\dT{}^*_P\mathbb{P}$, resp. $\dT{}^*_{P^*}\mathbb{P}^*$, the conormal space to $P$ in $\dT{}^*\mathbb{P}$, resp. to $P^*$ in $\dT{}^*\mathbb{P}^*$, and we will construct and denote by $\chi$ the homogeneous symplectic isomorphism between $\dT{}^*\mathbb{P}$ and $\dT{}^*\mathbb{P}^*$. 

For $\varepsilon\in\mathbb{Z}/2\mathbb{Z}$, we denote by $\mathbb{C}_{P}(\varepsilon)$ the following sheaves: for $\varepsilon=0$, we set
\eqn
\mathbb{C}_{P}(0):=\mathbb{C}_{P}
\eneqn
for $\varepsilon=1$, $\mathbb{C}_{P}(1)$ is the sheaf defined by the following exact sequence: 
\eq\label{def:proj_loc_constant_sheaf}
0\rightarrow \mathbb{C}_{P}(1) \rightarrow \eim{q}\mathbb{C}_{\widetilde{P}}\xrightarrow{tr}\mathbb{C}_{P}\rightarrow0
\eneq
where $q$ is the $2:1$ map from the universal cover $\widetilde{P}$ of $P$, to $P$ and $tr$ the integration morphism $tr:\eim{q}\mathbb{C}_{\widetilde{P}}\simeq\eim{q}\epb{q}\mathbb{C}_{P}\to\mathbb{C}_{P}$. 

Let an integer $p\in\mathbb{Z}$, $\varepsilon\in\mathbb{Z}/2\mathbb{Z}$, we define the sheaves of real analytic functions, hyperfunctions on $P$ resp. $P^*$ twisted by some power of the tautological line bundle, 
\eqn
\sha_{P}(\varepsilon,p)\eqdot\sha_{P}\tens_{\sho_\mathbb{P}}\sho_\mathbb{P}(p)\tens_{\mathbb{C}}\mathbb{C}_{P}(\varepsilon)
\eneqn
\vspace{-2em}
\eqn
\shb_P(\varepsilon,p):= \shb_P\tens_{\sha_{P}} \sha_{P}(\varepsilon,p)\simeq\rhom(\RD'_\mathbb{P}\C_P,\sho_\mathbb{P}(p))\tens\mathbb{C}_{P}(\varepsilon)
\eneqn

We define the sheaves of microfunctions on $P$ resp. $P^{*}$ twisted by some power of the tautological bundle, 
\eqn
\mathscr{C}_{P}(\varepsilon,p):=H^0(\mu\mathpzc{hom}(D'_\mathbb{P}\mathbb{C}_{P},\mathcal{O}_{\mathbb{P}}(p)))\tens\mathbb{C}_{P}(\varepsilon)
\eneqn
and similarly with $P^*$ instead of $P$. We notice that for $n$ odd, $D'_\mathbb{P}\mathbb{C}_{P}\simeq\mathbb{C}_{P}(0)=\mathbb{C}_{P}$, and for $n$ even $D'_\mathbb{P}\mathbb{C}_{P}\simeq\mathbb{C}_{P}(1)$.

For $X,Y$ either the manifold $\mathbb{P}$ or $\mathbb{P}^*$, for any two integers $p,q$, we note $\mathcal{O}_{X\times Y}(p,q)$ the line bundle on $X\times Y$ with homogenity $p$ in the $X$ variable and $q$ in the $Y$ variable. We set 
\eqn
\begin{array}{rrl}
	\Omega_{X\times Y/X}(p,q)&:=&\Omega_{X\times Y/X} \tens_{\mathcal{O}_{X\times Y}}\mathcal{O}_{X\times Y}(p,q)\\
	\she_X^{\mathbb{R}}(p,q)&:=&\muhom(\mathbb{C}_{\Delta_{X}},\Omega_{X\times X/X}(p,q))[d_X]
\end{array}
\eneqn
and we define accordingly $\she_X(p,q)$. Let us notice that $\she_X^{\mathbb{R}}(-p,p)$ is a sheaf of rings. 

Let $n$ be the dimension of $P$, (of course $n=d_{\mathbb{P}}$). For an integer $k$ and $\varepsilon\in\mathbb{Z}/2\mathbb{Z}$, we note 
\eqn
k^{*}:=-n-1-k
\eneqn 
\eqn
\varepsilon^{*}:=-n-1-\varepsilon\hspace{1ex}mod(2)
\eneqn

\subsection{Projective duality: geometry} 
\label{geometric_situation}

\subsubsection{Notations}\label{notation_tensored_microfunctions}

We refer to the notations of the sections~\ref{sec:notation_from_results}. We recall that we denote by 
\eqn
\text{$V$, $\mathbb{V}$, an $(n+1)$-dimensional  real and complex vector space, respectively,}\\
\text{$P$, $\mathbb{P}$, the $n-$dimensional real and complex projective space, respectively,}\\
\text{$S$, $\mathbb{S}$, the real and complex incidence hypersurface in $P\times P^*$, $\mathbb{P}\times\mathbb{P}^*$, respectively.}
\eneqn
When necessary, we will enforce the dimension by noting $\mathbb{P}_{n}$, resp. $\mathbb{P}^{*}_{n}$.

Let $X,Y$ be complex manifolds, we recall that we denote by $q_1$ and $q_2$ the respective projection of $X\times Y$ on each of its factor.

For $K\in\Derb(\mathbb{C}_{X\times Y})$, we recall that we defined the functor:
\eqn
&& \Phi_K \cl  \Derb(\mathbb{C}_X)\rightarrow\Derb(\mathbb{C}_Y) \\
&& F\mapsto Rq_{2 !}(K\tens q_{1}^{-1}F)\\
\eneqn

For an integer $k$ and $\varepsilon\in\mathbb{Z}/2\mathbb{Z}$, we note $k^{*}‹=-n-1-k$ and $\varepsilon^{*}=-n-1-\varepsilon \text{ mod}(2)$.  We refer to Section \ref{sec:notation_from_results} for the definition of the sheaves of twisted microfunctions  $\mathscr{C}_{P}(\varepsilon,k),\mathscr{C}_{P^{*}}(\varepsilon^*,k^*)$.

%

%

\subsubsection{Geometry of projective duality}

For a manifold $X$, we denote by $P^*X$ the projectivization of the cotangent bundle of $X$. The following results are well-known. However, we will give a proof of Prooposition \ref{dual_topo_iso} since it is more straighforward than the one usually found in the litterature.

\begin{proposition}\label{dual_topo_iso}
	There is an homogeneous complex symplectic isomorphism
	\eq\label{eq:dual_topo_cotangent_iso_t}
	 \dT{}^{*}\mathbb{P} \simeq \dT{}^{*}\mathbb{P}^{*}
	\eneq
	and a contact isomorphism
	\eq\label{eq:dual_topo_coproj_iso_t}
	P^*\mathbb{P}\simeq \mathbb{S}\simeq P^*\mathbb{P}^*
	\eneq
\end{proposition}

\begin{proof}
	We have the natural morphism
	\eqn
	\mathbb{V}\setminus\{0\} \xrightarrow{\rho} \mathbb{P}
	\eneqn
	According to \ref{diag:microlocal1}, this morphism, after removing the zero section,  induces the following diagram
	\eqn
	&&\xymatrix{
	T{}^*(\mathbb{V}\setminus\{0\}) \ar[dr]&\mathbb{V}\setminus\{0\} \times_{\mathbb{P}} T{}^*\mathbb{P}\ar[l]_-{^t\rho'}\ar[r]\ar[d]&T{}^*\mathbb{P}\ar[d]\\
	&\mathbb{V}\setminus\{0\}\ar[r]&\mathbb{P}
	}
	\eneqn
	We notice that $^t\rho'$ is an immersion. Let us denote by $\mathbb{H}$, $\mathbb{H}^*$, the incidence hypersurfaces: 
	\eqn
	\mathbb{H}=\{(\xi,x)\in\mathbb{V}^*\times(\mathbb{V}\setminus\{0\});\langle\xi,x\rangle=0\}\\
	\mathbb{H}^*=\{(x,\xi)\in\mathbb{V}\times(\mathbb{V}^*\setminus\{0\});\langle x,\xi\rangle=0\}
	\eneqn
	Noticing that for $x\in\mathbb{V}\setminus\{0\}$, $\rho$ is constant along the fiber above $\rho(x)$, we see that $^t\rho'$ is an immersion into the incidence hypersurface $\mathbb{H}$. Besides, $^t\rho'$ is a morphism of fibered space and so, by a dimensional argument, we conclude that this immersion is also onto.
	
	Removing the zero sections, we get the diagram 
	\eq\label{dia:incidence_map_duality}
	&&\xymatrix{
	T{}^*(\mathbb{V}\setminus\{0\}) \ar[d]_-{\simeq} & \mathbb{H} \ar@{_{(}->}[l] \ar[d]_-{\simeq}  &\mathbb{V}\setminus\{0\} \times_{\mathbb{P}} \dT{}^*\mathbb{P}\ar[l]\ar[r]\ar[d]_-{\simeq} & \dT{}^*\mathbb{P}\\
	T{}^*(\mathbb{V}^*\setminus\{0\}) & \mathbb{H}^* \ar@{_{(}->}[l] & (\mathbb{V}^*\setminus\{0\} )\times_{\mathbb{P}^*} \dT{}^*\mathbb{P}^*  \ar[l]\ar[r]& \dT{}^*\mathbb{P}^*
	}
	\eneq
	where the isomophism between $\mathbb{H}$ and $\mathbb{H}^*$ follows from the following symplectic isomorphism:
	\eqn
	\dT{}^*(\mathbb{V}\setminus\{0\}) \simeq \dT{}^*(\mathbb{V}^*\setminus\{0\}) \\
	(x,\xi) \mapsto (\xi,-x)
	\eneqn
	Now, taking the quotient by the action of $\mathbb{C}^*$ on both sides of the isomorphism between $(\mathbb{V}\setminus\{0\} )\times_{\mathbb{P}} \dT{}^*\mathbb{P}$ and $(\mathbb{V}^*\setminus\{0\} )\times_{\mathbb{P}^*} \dT{}^*\mathbb{P}^*$, we get the isomorphism: 
	\eqn
	 \dT{}^*\mathbb{P} \simeq (\mathbb{V}^*\setminus\{0\} )\times_{\mathbb{P}^*} P^*\mathbb{P}^* \simeq \dT{}^*\mathbb{P}^*
	\eneqn
	
	This gives (\ref{eq:dual_topo_cotangent_iso_t}).
	
	Besides, passing to the quotient by the action of $\mathbb{C}^*\times\mathbb{C}^*$ on the two central columns of diagram (\ref{dia:incidence_map_duality}), we get (\ref{eq:dual_topo_coproj_iso_t}).
\end{proof}
\begin{proposition}\label{projective_duality_iso}
	Consider the double fibrations
	\eq\label{dia:projective_duality_iso}
	&&\xymatrix{
		&\dT{}^{*}_{\mathbb{S}}(\mathbb{P}\times \mathbb{P}^{*})\ar[ld]_\sim^-{p_1}\ar[rd]^\sim_-{p_{2}^{a}}&\\
		\dT{}^{*}\mathbb{P}\ar[rr]^-\sim_-\chi&  & \dT{}^{*}\mathbb{P}^{*} 
	}
	\eneq	
	
	Then, $p_{1}$ and $p_{2}^{a}$ are isomorphisms and $\chi=p_{2}^{a} \circ \opb{p_{1}}$ is a homogeneous symplectic isomorphism.
\end{proposition}

Now, we are going to prove the following

\begin{proposition}\label{rc_dual_topo_iso}
	The diagram \ref{dia:projective_duality_iso} induces
	\eqn
	&&\xymatrix{
		& \dT{}^{*}_{\mathbb{S}}(\mathbb{P}\times \mathbb{P}^{*})\cap(\dT{}^{*}_{P}\mathbb{P}\times\dT{}^{*}_{P^{*}}\mathbb{P}^{*})\ar[ld]_\sim^-{p_1}\ar[rd]^\sim_-{p_{2}^{a}}&\\
		\dT{}^{*}_{P}\mathbb{P}\ar[rr]^-\sim_-\chi&  & \dT{}^{*}_{P^{*}}\mathbb{P}^{*} 
	}
	\eneqn		
\end{proposition}

\subsection{Projective duality for microdifferential operators}
\label{presentation_preparatory_theorem}

Let $k,k'$ be integers and $\varepsilon\in\mathbb{Z}/2\mathbb{Z}$. We follow the notations of the sections \ref{sec:notation_from_results} and \ref{sec:reminder_algebraic_analysis}. We define similarly a twisted version of $\shb_{\mathbb{S}\vert \mathbb{P}\times \mathbb{P}^*}^{(n,0)}$ and $\shc_{\mathbb{S}\vert \mathbb{P}\times \mathbb{P}^*}^{(n,0)}$.

We set 
\eqn
\shb_{\mathbb{S}}^{(n,0)}(k,k'):=\opb{q_{2}}\mathcal{O}_{\mathbb{P}^*}(k')\tens_{\opb{q_2}\mathcal{O}_{\mathbb{P}^*}}\shb_{\mathbb{S}\vert \mathbb{P}\times \mathbb{P}^*}\tens_{\opb{q_1}\mathcal{O}_\mathbb{P}}\opb{q_{1}}(\mathcal{O}_{\mathbb{P}}(k)\tens_{\mathcal{O}_\mathbb{P}}\Omega_\mathbb{P})
\eneqn

and the $(\she_{\mathbb{P}}(-k,k),\she_{\mathbb{P}^*}(-k^*,k^*))$-module 
\eqn
\shc_{\mathbb{S}\vert \mathbb{P}\times \mathbb{P}^*}^{(n,0)}(k,k'):=\she \shb_{\mathbb{S}\vert \mathbb{P}\times \mathbb{P}^*}^{(n,0)}(k,k')
\eneqn

We notice that $\she_{\mathbb{P}}(-k,k)$ is nothing but $\mathcal{O}_{\mathbb{P}}(-k)\mathcal{D}\tens{_{\opb{\pi{_{\mathbb{P}}}}\mathcal{D}_\mathbb{P}}}\she_{\mathbb{P}}\tens{_{\opb{\pi{_{\mathbb{P}}}}\mathcal{D}_\mathbb{P}}}\mathcal{D}\mathcal{O}_{\mathbb{P}}(k)$. According to the diagram \ref{dia:projective_duality_iso}, we denoted by $\chi$ the homogeneous symplectic isomorphism 
\eqn
\chi:=p_{2}^{a}\vert_{\dT{}^{*}_{\mathbb{S}}(\mathbb{P}\times \mathbb{P}^{*})}\circ\opb{p_{1}\vert_{\dT{}^{*}_{\mathbb{S}}(\mathbb{P}\times \mathbb{P}^{*})}}
\eneqn

We have: 
\begin{theorem}[{\cite[p.~469]{DS96}}]\label{th:projective_non_degenerate_section}
	Assume $-n-1<k<0$. There exists a section $s$ of $\muhom(\mathbb{C}_{\mathbb{S}}[-1],\Omega_{\mathbb{P}\times \mathbb{P}^*/\mathbb{P}^*}(-k,k^*))$, non-degenerate on $\dT{}^{*}_{\mathbb{S}}(\mathbb{P}\times\mathbb{P}^*)$. 
\end{theorem}
\begin{proof}
	 From the exact sequence:
	 \eq\label{support_distinguished_triangle}
	 0 \longrightarrow \mathbb{C}_{(\mathbb{P}\times\mathbb{P}^{*})\setminus \mathbb{S}}\longrightarrow \mathbb{C}_{\mathbb{P}\times\mathbb{P}^{*}} \longrightarrow \mathbb{C}_{\mathbb{S}} \longrightarrow 0
	 \eneq
	 we get the natural morphism
	 \eqn
	 \begin{array}{rcl}
	 	\rsect((\mathbb{P}\times\mathbb{P}^{*})\setminus \mathbb{S};\Omega_{\mathbb{P}\times \mathbb{P}^*/\mathbb{P}^*}(-k,k^*))\ & \to & \rsect_{\mathbb{S}}(\mathbb{P}\times\mathbb{P}^{*};\Omega_{\mathbb{P}\times \mathbb{P}^*/\mathbb{P}^*}(-k,k^*))[1] \\
	 	& \simeq & \rsect(\mathbb{P}\times\mathbb{P}^{*};\rhom(\mathbb{C}_{\mathbb{S}};\Omega_{\mathbb{P}\times \mathbb{P}^*/\mathbb{P}^*}(-k,k^*))) [1]\\
	 	& \simeq & \rsect(T^*{}(\mathbb{P}\times\mathbb{P}^*);\muhom(\mathbb{C}_{\mathbb{S}};\Omega_{\mathbb{P}\times \mathbb{P}^*/\mathbb{P}^*}(-k,k^*))) [1]\\
	 	& \to &
	 	\rsect(\dT{}^*(\mathbb{P}\times\mathbb{P}^*);\muhom(\mathbb{C}_{\mathbb{S}}[-1];\Omega_{\mathbb{P}\times \mathbb{P}^*/\mathbb{P}^*}(-k,k^*)))\\
	 \end{array}
	 \eneqn
	 
	 Let $z=(z_0,... ,z_n)$ be a system of homogeneous coordinates on $\mathbb{P}$ and $\zeta= (\zeta_0,... ,\zeta_n)$ the dual coordinates on $\mathbb{P}^*$. As explained in \cite{DS96}, a non-degenerate section is provided by the Leray section, defined for $(z,\xi)\in(\mathbb{P}\times\mathbb{P}^*)\setminus\mathbb{S}$ by
	\eq\label{def:projective_non_degenerate_section}
	s(z,\zeta)=\frac{\omega'(z)}{\langle z,\zeta\rangle^{n+1+k}}
	\eneq
	where $\omega'(z)$ is the Leray form $\omega^{\prime }(z) = \sum_{k=0}^{n}(-1)^{k}z_{k} dz_{0} {\small %
	\wedge }\ldots\wedge dz_{k-1}\wedge dz_{k+1}\wedge\ldots{\small \wedge }dz_{n}$, Leray~\cite{L59}.
\end{proof}

Let $s$ be a section of $H^1(\muhom(\mathbb{C}_{\mathbb{S}}[-1],\Omega_{\mathbb{P}\times \mathbb{P}^*/\mathbb{P}^*}(-k,k^*)))$, non-degenerate on $\dT{}^{*}_{\mathbb{S}}(\mathbb{P}\times\mathbb{P}^*)$.
\begin{theorem}\label{th:projective_operator_main_theorem}
	Assume $-n-1<k<0$. Then, we have an isomorphism in $\Derb(\C_{\dT{}^*\mathbb{P}})$
	\eqn
	\underline{\Phi}_{\mathbb{S}}^{\mu}(\she_{\mathbb{P}}(-k,k)\vert_{\dT{}^*\mathbb{P}}) \simeq \she_{\mathbb{P}^*}(-k^*,k^*)\vert_{\dT{}^*\mathbb{P}^*}\\
	\oim{\chi}\she_{\mathbb{P}}(-k,k)\vert_{\dT{}^*\mathbb{P}} \simeq \she_{\mathbb{P}^*}(-k^*,k^*)\vert_{\dT{}^*\mathbb{P}^*}
	\eneqn
\end{theorem}
\begin{proof}
	Let $\mathcal{F}$, $\mathcal{G}$ be line bundles on $\mathbb{P}$, and $\mathbb{P}^*$ respectively. We know from \cite{SKK73} that a global non-degenerate section $s\in\sect(\dT{}^*\mathbb{P}\times\dT{}^*\mathbb{P}^*;\shc_{\mathbb{S}\vert \mathbb{P}\times \mathbb{P}^*}^{(n,0)}\tens_{\opb{p_1}\she_{\mathbb{P}}}\she\mathcal{F}\tens_{\opb{p_2}\she_{\mathbb{P}^*}}\mathcal{G}^{\tens -1}\she)$ induces an isomorphism  of $\she$-modules
	\eqn
	\underline{\Phi}_{\mathbb{S}}^{\mu}(\she\mathcal{F}\vert_{\dT{}^*\mathbb{P}}) \simeq \she\mathcal{G}\vert_{\dT{}^*\mathbb{P}^*}
	\eneqn
	Now, let us set $\mathcal{F}=\mathcal{O}_{\mathbb{P}}(k)$, $\mathcal{G}=\mathcal{O}_{\mathbb{P}^{*}}(k^{*})$. Then, \ref{th:projective_non_degenerate_section} provides such a non-degenerate section in $\sect(\dT{}^*\mathbb{P}\times\dT{}^*\mathbb{P}^*;\shc_{\mathbb{S}\vert \mathbb{P}\times \mathbb{P}^*}^{(n,0)}\tens_{\opb{p_1}\she_{\mathbb{P}}}\she\mathcal{O}_{\mathbb{P}}(k)\tens_{\opb{p_2}\she_{\mathbb{P}^*}}\mathcal{O}_{\mathbb{P}^{*}}(k^{*})^{\tens -1}\she)$. So that, we have an isomorphism
	\eqn
	\underline{\Phi}_{\mathbb{S}}^{\mu}(\she_{\mathbb{P}}(-k,k)\vert_{\dT{}^*\mathbb{P}}) \simeq \she_{\mathbb{P}^*}(-k^*,k^*)\vert_{\dT{}^*\mathbb{P}^*}
	\eneqn
	On the other hand $s$ is a non-degenerate section of $\shc_{\mathbb{S}\vert \mathbb{P}\times \mathbb{P}^*}^{(n,0)}(-k,k^*)$, hence we can apply Theorem \ref{th:quantized_iso}. Let us denote by 
	\eqn
	\she_{\mathbb{P}\times \mathbb{P}^*}(k,k^*):=\she_{\mathbb{P}}(-k,k)\etens_{\opb{\pi}\mathcal{O}_{\mathbb{P}\times\mathbb{P^*}}}\she_{\mathbb{P}^*}(-k^*,k^*)
	\eneqn
	Theorem \ref{th:quantized_iso} gives the following isomorphisms
	\eqn
	\she_{\mathbb{P}}(-k,k)\vert_{\dT{}^*\mathbb{P}} \simeq \oim{p_1}(\she_{\mathbb{P}\times \mathbb{P}^*}(k,k^*).s)\vert_{\dT{}^*\mathbb{P}}\\
	\oim{p_2}(\she_{\mathbb{P}\times \mathbb{P}^*}(k,k^*).s)\vert_{\dT{}^*\mathbb{P}^*} \simeq
	\she_{\mathbb{P}^*}(-k^*,k^*)\vert_{\dT{}^*\mathbb{P}^*}
	\eneqn
	And so
	\eqn
	\oim{\chi}\she_{\mathbb{P}}(-k,k)\vert_{\dT{}^*\mathbb{P}} \simeq \she_{\mathbb{P}^*}(-k^*,k^*)\vert_{\dT{}^*\mathbb{P}^*}
	\eneqn
\end{proof}

\subsection{Projective duality for microfunctions }
\label{presentation_main_theorem}

In the following, we will denote by $K$ the object $\mathbb{C}_{\mathbb{S}}[n-1]$. In order to prove Proposition \ref{microfunction_iso_theorem}, we will need to compute $\Phi_{K}(\mathbb{C}_{P}(1))$, which is done in \cite{DS96}:

\begin{lemma}[{\cite{DS96}}]\label{lem:constant_sheaf_correspondance}
	We have
	\eqn
	\Phi_K(\mathbb{C}_{P}(1))\simeq\left\{
	\begin{array}{ll}
	\mathbb{C}_{P^{*}}(1)\text{ , for $n$ odd}\\
	\mathbb{C}_{\mathbb{P}^{*}\setminus P^{*}}[1]\text{ , for $n$ even}\\
	\end{array}
	\right.
	\eneqn
	and
	\eqn
	H^j(\Phi_{K}(\mathbb{C}_{P}(0)))\simeq\left\{
	\begin{array}{ll}
	\mathbb{C}_{\mathbb{P}^*}\text{ , for $j=n-1$}\\
	\mathbb{C}_{\mathbb{P}^*\setminus P^*}\text{, for $j=-1$ and $n$ odd}\\
	\mathbb{C}_{P^*}(1)\text{, for $j=0$ and $n$ even}\\
	0\text{ in any other case}\\
	\end{array}
	\right.
	\eneqn
\end{lemma}
We are in a proposition to prove:

\begin{theorem}\label{microfunction_iso_theorem}
	Assume $-n-1<k<0$. Recall that any section $s\in \sect(\mathbb{P}\times\mathbb{P}^*;\shb_{\mathbb{S}}^{(n,0)}(-k,k^*))$, defines a morphism in $\Derb(\C_{\dT\mathbb{P}})$
	\eq\label{projective_twisted_microfunction_iso}
	\chi_{*}\mathscr{C}_{P}(\varepsilon,k)\vert_{\dT{}^*_P\mathbb{P}}\to\mathscr{C}_{P^{*}}(\varepsilon^{*},k^{*})\vert_{\dT{}^*_{P^*}\mathbb{P}^*}
	\eneq
	
	Assume $s$ is non-degenerate on $\dT{}^*_{\mathbb{S}}(\mathbb{P}\times\mathbb{P}^*)$. Then (\ref{projective_twisted_microfunction_iso}) is an isomorphism. Moreover, there exists such a non-degenerate section.
\end{theorem}

\begin{remark}
	\bnum
	\item This is a refinement of a general theorem of \cite{SKK73} and is a microlocal version of Theorem 5.17 in \cite{DS96}.
	\item 
	The classical Radon transform deals with the case where $k=-n$, $k^{*}=-1$.
	\enum

\end{remark}

\begin{proof}
	We will deal with the case $\varepsilon=1$ and $n$ even, the complementary cases being proven the same way. Let us apply Theorem $\ref{main_theorem_contact_muhom}$ in the following particular case
	
	- $U=\dT{}^{*}\mathbb{P}$, $V=\dT{}^{*}\mathbb{P}^*$, $\Lambda=\dT{}^{*}_{\mathbb{S}}(\mathbb{P}\times\mathbb{P}^{*})$.
	
	- $K$ is $\mathbb{C}_{\mathbb{S}}[n-1]$.
	
	- $F_{1}=\mathbb{C}_{P}(1)$ and $F_{2}=\mathcal{O}_{\mathbb{P}}(k)$
	
	$K$ verifies conditions (i),(ii),(iii) of Theorem $\ref{main_theorem_contact_muhom}$
	
	(i) is fulfilled as the constant sheaf on a closed submanifold of a manifold is cohomologically constructible.
	
	(ii) is fulfilled since $SS(\mathbb{C}_{\mathbb{S}})$ is nothing but $T^{*}_{\mathbb{S}}(\mathbb{P}\times\mathbb{P}^{*})$.
	
	(iii) $\mathbb{C}_{T^{*}_{\mathbb{S}}(\mathbb{P}\times\mathbb{P}^{*})}\longrightarrow \mu\mathpzc{hom}(\mathbb{C}_{\mathbb{S}},\mathbb{C}_{\mathbb{S}})$ is an isomorphism on $T^{*}_{\mathbb{S}}(\mathbb{P}\times\mathbb{P}^{*})$ (this follows from the fact that for a closed submanifold $Z$ of a manifold $X$, $\mu_{Z}(\mathbb{C}_Z)\isoto\mathbb{C}_{T^*_ZX}$, see \cite[Prop.~4.4.3]{KS90}).

	 By a fundamental result in~\cite[Th 5.17]{DS96}, we know that for $-n-1<k<0$, a section $s\in \sect(\mathbb{P}\times\mathbb{P}^*;\shb_{\mathbb{S}}^{(n,0)}(-k,k^*))$, non-degenerate on $\dT{}^*_{\mathbb{S}}(\mathbb{P}\times\mathbb{P}^*)$, induces an isomorphism
	\eqn
	\Phi_K(\mathcal{O}_{\mathbb{P}}(k))\simeq \mathcal{O}_{\mathbb{P}^{*}}(k^{*})
	\eneqn
	Formula (\ref{def:projective_non_degenerate_section}) provides an example of such a non-degenerate section. Hence, applying Lemma \ref{lem:constant_sheaf_correspondance}, Theorem $\ref{main_theorem_contact_muhom}$ gives:
	\eqn
	\begin{split}
	\oim{\chi}\mu\mathpzc{hom}(\mathbb{C}_{P}(1),\mathcal{O}_{\mathbb{P}}(k))\vert_{\dT{}^*\mathbb{P}} & \simeq\mu\mathpzc{hom}(\mathbb{C}_{\mathbb{P}^{*}\setminus P^{*}}[1],\mathcal{O}_{\mathbb{P}^{*}}(k^{*}))\vert_{\dT{}^*\mathbb{P}^{*}}\\
	\end{split}
	\eneqn

	We have the exact sequence:
	\eq\label{support_distinguished_triangle}
	0 \longrightarrow \mathbb{C}_{\mathbb{P}^{*}\setminus P^*} \longrightarrow \mathbb{C}_{\mathbb{P}^{*}} \longrightarrow \mathbb{C}_{P^*} \longrightarrow 0
	\eneq
	Now, for any $F\in\Derb(\mathbb{C}_{\mathbb{P}^*})$, we have 
	\eqn
	supp(\muhom(\mathbb{C}_{\mathbb{P}^*},F)\vert_{\dT{}^*\mathbb{P}^{*}})\subset (SS(\mathbb{C}_{\mathbb{P}^*}) \cap \dT{}^*\mathbb{P}^{*}) \cap SS(F) = \emptyset
	\eneqn
	and hence, 
	\eqn
	\muhom(\mathbb{C}_{\mathbb{P}^*},F)\vert_{\dT{}^*\mathbb{P}^{*}}\simeq 0
	\eneqn

	Applying the $\muhom$ functor to $\ref{support_distinguished_triangle}$, we get
	\eqn
	\muhom(\mathbb{C}_{\mathbb{P}^{*}\setminus P^{*}},F)\vert_{\dT{}^*\mathbb{P}^{*}}[-1]\simeq \muhom(\mathbb{C}_{P^{*}},F)\vert_{\dT{}^*\mathbb{P}^{*}}
	\eneqn
	
	Hence, we have proved in particular that
	\eqn
	\begin{split}
	\oim{\chi}\mu\mathpzc{hom}(\mathbb{C}_{P}(1),\mathcal{O}_{\mathbb{P}}(k))\vert_{\dT{}^*_P\mathbb{P}} & \simeq \mu\mathpzc{hom}(\mathbb{C}_{P^{*}},\mathcal{O}_{\mathbb{P}^{*}}(k^{*}))\vert_{\dT{}^*_{P^*}\mathbb{P}^{*}}\\
	\end{split}
	\eneqn
\end{proof}

\subsection{Main results}

We follow the notations of Section~\ref{sec:notation_from_results} and Section~\ref{sec:reminder_algebraic_analysis}. 

Let us consider the situation (\ref{dia:projective_duality_iso}), where we denoted by $\chi$ the homogeneous symplectic isomorphism between $\dT{}^*\mathbb{P}$ and $\dT{}^*\mathbb{P}^*$ through $\dT{}^*_\mathbb{S}(\mathbb{P}\times\mathbb{P}^*)$. We set 
\eqn
L:=\mathbb{C}_{\mathbb{S}}[-1]
\eneqn
Then $L$ is a perverse sheaf satisfying
\eq\label{eq:2}
&&(\opb{p_1}(\dT{}^*\mathbb{P})\cup\opb{{p_2^a}}(\dT{}^*\mathbb{P}^*))\cap\SSi(L)=\dT{}^*_\mathbb{S}(\mathbb{P}\times\mathbb{P}^*)
\eneq

Recall Theorem \ref{th:projective_non_degenerate_section}, and let $s$ be a section of $\muhom(\mathbb{C}_{\mathbb{S}}[-1],\Omega_{\mathbb{P}\times \mathbb{P}^*/\mathbb{P}^*}(-k,k^*))$, non-degenerate on $\dT{}^{*}_{\mathbb{S}}(\mathbb{P}\times\mathbb{P}^*)$. We are in situation to apply Theorem \ref{th:qct_main_theorem}.
\begin{theorem}\label{th:projective_main_theorem}
	Let $G\in\Derb(\C_{\mathbb{P}^*})$, $k$ an integer. Assume $-n-1<k<0$. Then, we have an isomorphism in $\Derb(\C_{\dT{}^*\mathbb{P}})$:
	\eq\label{eq:qct_projective_main_theorem}
	&&\opb{\chi}\muhom(G,\sho_{\mathbb{P}^*}(k^*))\isoto\muhom(\Phi_{\mathbb{C}_{\mathbb{S}}[n-1]}(G),\sho_{\mathbb{P}}(k))
	\eneq
	This isomorphism is compatible with the action of $\she_{\mathbb{P}^*}(-k^*,k^*)$ and $\she_{\mathbb{P}}(-k,k)$ on the left and right side of (\ref{eq:qct_projective_main_theorem}) respectively.
\end{theorem}
\begin{proof}
	The isomorphism is directly provided by Theorem \ref{th:qct_main_theorem} in the situation where, using the notation inthere, $U=\dT{}^*\mathbb{P}$, $V=\dT{}^*\mathbb{P}^*$ and $\Lambda=\dT{}^*_\mathbb{S}(\mathbb{P}\times\mathbb{P}^*)$ and where we twist by homogenous line bundles of $\mathbb{P}$, $\mathbb{P}^*$ as explained below.

	 Let us adapt (\ref{eq:qct_main_theorem}) by taking into account the twist by homogeneous line bundles. We follow the exact same reasoning than sections of \ref{sec:kernels_complex_manifolds} and \ref{sec:qct_main_theorem}. 
	 
	 We have the natural  morphism
	 	\eqn
	 	&&\Omega_{\mathbb{P}^*\times \mathbb{P}/\mathbb{P}}(-k^*,k)\conv\sho_{\mathbb{P}^*}(k^*)\,[n]\to\sho_{\mathbb{P}}(k).
	 	\eneqn

	Indeed, we have
	 	\eqn
	 	\Omega_{\mathbb{P}^*\times \mathbb{P}/\mathbb{P}}(-k^*,k)\conv\sho_{\mathbb{P}^*}(k^*)\,[n]&=&
	 	\reim{q_1}(\mathcal{O}_{\mathbb{P}^*\times \mathbb{P}}(-k^*,k)\tens_{\opb{q_2}\mathcal{O}_{\mathbb{P}^*}}\opb{q_2}\Omega_{\mathbb{P}^*}\tens\opb{q_2}\sho_{\mathbb{P}^*}(k^*)[n])\\
	 	&\to&\reim{q_1}(\mathcal{O}_{\mathbb{P}^*\times \mathbb{P}}(k,0)\tens_{\opb{q_2}\mathcal{O}_{\mathbb{P}^*}}\opb{q_2}\Omega_{\mathbb{P}^*})[n]\to[\int]\sho_{\mathbb{P}}(k)
	 	\eneqn


	 	Given this morphism and considering $\shl\eqdot\muhom(\mathbb{C}_{\mathbb{S}}[-1],\Omega_{\mathbb{P}^*\times\mathbb{P}/\mathbb{P}}(-k^*,k))$, we mimic the proof of Theorem \ref{th:KS14b} so that for a section $s$ of $\shl$ on $\dT{}^*\mathbb{P}\times\dT{}^*\mathbb{P}^{*a}$ and for $P\in\sect(\dT{}^*\mathbb{P};\she_{\mathbb{P}}(-k,k))$  and $Q\in\sect(\dT{}^*\mathbb{P}^{*};\she_{\mathbb{P}^*}(-k^*,k^*))$ satisfying $P\cdot s=s\cdot Q$, the diagram below commutes:
 	\eq\label{diag:Ps=sQ}
 	\eneq
 	\eqn
 	&&\xymatrix{
 	\mathbb{C}_{\mathbb{S}}\circ\muhom(G,\sho_{\mathbb{P}^*}(k^*))\vert_{\dT{}^*_{P^*}\mathbb{P}^*}\ar[rr]_-{\alpha(s)}\ar[d]_-{\Phi_{\mathbb{C}_{\mathbb{S}}}(\alpha(Q))}&&\muhom(\mathbb{C}_{\mathbb{S}}[n-1]\conv G,\sho_{\mathbb{P}}(k))\vert_{\dT{}^*_P\mathbb{P}}\ar[d]_-{\alpha(P)}\\
 	\mathbb{C}_{\mathbb{S}}\circ\muhom(G,\sho_{\mathbb{P}^*}(k^*))\vert_{\dT{}^*_{P^*}\mathbb{P}^*}\ar[rr]_-{\alpha(s)}&&\muhom(\mathbb{C}_{\mathbb{S}}[n-1]\conv G,\sho_{\mathbb{P}}(k))\vert_{\dT{}^*_P\mathbb{P}}.
 	}
 	\eneqn
 	
 	From there, given a non-degenerate section of $\shl$ on $\dT{}^{*}_S(\mathbb{P}\times\mathbb{P}^*)$, Theorem \ref{th:qct_main_theorem} gives the compatible action of micro-differential operators on each side of the isomorphism (\ref{eq:qct_projective_main_theorem})
 	\eqn
 	&&\opb{\chi}\muhom(G,\sho_{\mathbb{P}^*}(k^*))\vert_{\dT{}^*_{P^*}\mathbb{P}^*}\isoto\muhom(\Phi_{\mathbb{C}_{\mathbb{S}}[n-1]}(G),\sho_{\mathbb{P}}(k))\vert_{\dT{}^*_P\mathbb{P}}
 	\eneqn
	 
	It remains to exhibit a non-degenerate section so that, for $P\in\sect(\dT{}^*\mathbb{P};\she_{\mathbb{P}}(-k,k))$, there is $Q\in\sect(\dT{}^*\mathbb{P}^{*};\she_{\mathbb{P}^*}(-k^*,k^*))$ such that $P\cdot s=s\cdot Q$. Precisely, $s$ is given by Proposition \ref{th:projective_non_degenerate_section}.
\end{proof}

Specializing the above proposition, we get
\begin{corollary}\label{coro:microfunction_duality_iso}
	Let $\varepsilon\in\mathbb{Z}/2\mathbb{Z}$. In the situation of Proposition \ref{th:projective_main_theorem}, we have the isomorphism, compatible with the respective action of $\opb{p_{1}}\she_{\mathbb{P}}(-k,k)$ and $\opb{{p_{2}^a}}\she_{\mathbb{P}^*}(-k^*,k^*)$
	\eqn
	\chi_{*}\mathscr{C}_{P}(\varepsilon,k)\vert_{\dT{}^*_P\mathbb{P}}\simeq\mathscr{C}_{P^{*}}(\varepsilon^{*},k^{*})\vert_{\dT{}^*_{P^*}\mathbb{P}^*}
	\eneqn
	 
\end{corollary}
\begin{proof}
	This is an immediate consequence of Proposition \ref{th:projective_main_theorem}, where we consider the special case $G=\mathbb{C}_{P^*}(\varepsilon^*)$. Indeed, we have, from Lemma \ref{lem:constant_sheaf_correspondance}, the isomorphism in $\Derb(\mathbb{C}_{\mathbb{P}^*};\dT{}^*\mathbb{P}^*)$
	\eqn
	\mathbb{C}_{\mathbb{S}}[n-1]\conv \mathbb{C}_{P^*}(\varepsilon^*)\simeq\mathbb{C}_{P}(\varepsilon)
	\eneqn 
\end{proof}


We can state now

\begin{corollary}
	 Let $k$ be an integer. Let $\mathcal{N}$ be a coherent $\she_{\mathbb{P}}(-k,k)$-module and $F\in\Derb(\mathbb{P})$.  Assume $-n-1<k<0$. Then, we have an isomorphism in $\Derb(\C_{\dT{}^*\mathbb{P}})$
	\eqn\hspace{-19em}
	\oim{\chi}\rhom{_{\she_{\mathbb{P}}(-k,k)}}(\mathcal{N},\muhom(F,\mathcal{O}_\mathbb{P}(k)))\simeq
	\eneqn
	\vspace{-2em}
	\eqn
	\hspace{15em}\rhom{_{\she_{\mathbb{P}^*}(-k^*,k^*)}}(\underline{\Phi}_{\mathbb{S}}^{\mu}(\mathcal{N}),\muhom((\Phi_{\mathbb{C}_{\mathbb{S}}[n-1]}F,\mathcal{O}_{\mathbb{P}^*}(k^*)))
	\eneqn
\end{corollary}
\begin{proof}
	It suffices to prove this statement for finite free $\she_{\mathbb{P}}(-k,k)$-modules, which in turn can be reduced to the case where $\mathcal{N}=\she_{\mathbb{P}}(-k,k)$. 
	By Theorem \ref{th:projective_operator_main_theorem}, we have
	\eqn
	\underline{\Phi}_{\mathbb{S}}^{\mu}(\she_{\mathbb{P}}(-k,k)\vert_{\dT{}^*\mathbb{P}}) \simeq \she_{\mathbb{P}^*}(-k^*,k^*)\vert_{\dT{}^*\mathbb{P}^*}
	\eneqn

	Then, by applying Proposition \ref{th:projective_main_theorem}, we have
	\eqn
	\oim{\chi}\muhom(F,\mathcal{O}_\mathbb{P}(k))\vert_{\dT{}^*\mathbb{P}^*} \simeq
	\rhom{_{\she_{\mathbb{P}^*}(-k^*,k^*)}}(\she_{\mathbb{P}^*}(-k^*,k^*),\muhom((\Phi_{\mathbb{C}_{\mathbb{S}}[n-1]}F,\mathcal{O}_{\mathbb{P}^*}(k^*)))\vert_{\dT{}^*\mathbb{P}^*}
	\eneqn
	which proves the corollary.
\end{proof}

\printbibliography[heading=bibintoc,title={References}]

\end{document}